\documentclass[11pt]{article}
\newcommand		{\comment}[1]		{}

\usepackage{amsmath,amsthm}
\usepackage{mathtools}

\usepackage{xspace}
\usepackage{enumitem}

\usepackage[left=2.54cm,right=2.54cm,top=2.54cm,bottom=2.54cm]{geometry}
\usepackage[markcase=ignoreuppercase]{scrlayer-scrpage}
\usepackage{xspace}
\usepackage{enumitem}

\usepackage{avant} 				
\usepackage[sc,osf]{mathpazo} 			
	\AtBeginDocument{
		\DeclareSymbolFont{AMSb}{U}{msb}{m}{n}
		\DeclareSymbolFontAlphabet{\mathbb}{AMSb}
	}
\usepackage{extpfeil} 
\usepackage{mathabx}  
\usepackage{mathrsfs}	
\usepackage{fancyvrb}	

\usepackage{extpfeil} 
\usepackage{mathabx}  
\usepackage{mathrsfs}	
\usepackage{fancyvrb}	

\usepackage{graphicx}
\usepackage{float}
\usepackage{subcaption}
\usepackage{caption}
\newcommand{\mockalph}[1]{\!}

\usepackage[utf8]{inputenc}
\usepackage[T1]{fontenc}

\usepackage[hyphens]{url}		%
\usepackage{hyperref}	
\usepackage{cleveref} 	
\newcommand		{\myred}		{BrickRed}
\hypersetup{
			colorlinks		= true, 	
			urlcolor		= \myred, 
			linkcolor		= \myred, 	
			citecolor		= PineGreen 	
}
\newcommand		{\hyref}[1]		{\hyperref[#1]{\ref*{#1}}}

\newif\ifdebug
\ifdebug\usepackage{lineno}\linenumbers\else\fi

\usepackage{etoolbox}
\PassOptionsToPackage{dvipsnames,svgnames,x11names,table}{xcolor}
\newcommand		{\defd}[1]	{\textcolor{RoyalBlue}{\textbf{\textit{#1}}}}
\newcommand		{\defm}[1]	{\textcolor{RoyalBlue}{#1}}

\newcommand		{\NB}		{\nd\textcolor{Red}{N.B.}}

\usepackage[all]{xypic}\xyoption{rotate}
	\newdir{ >}{{}*!/-7pt/\dir{>}}
	\newdir{i}{{}*!/-7pt/\dir{(}	}
\usepackage{tikz}
\usetikzlibrary{matrix,fit,backgrounds,calc,arrows} 
\tikzstyle{image}=[rectangle,fill=Red!20,inner sep=-2pt]
\tikzstyle{nonzero}=[rectangle,fill=Navy!20,inner sep=0pt]
\tikzstyle{nonzerosm}=[rectangle,fill=Navy!20,inner sep=-2pt]

\usepackage{thmtools}
\makeatletter
\let\c@figure\c@table
\let\c@equation\c@table
\makeatother

\numberwithin{table}{section}
\numberwithin{figure}{section}
\newtheorem{abcthm}[table]{Theorem}

\newtheorem{theorem}[table]{Theorem}

\newtheorem{proposition}[table]{Proposition}

\newtheorem{corollary}[table]{Corollary}

\newtheorem{lemma}[table]{Lemma}

\newtheorem{claim}[table]{Claim}

\theoremstyle{definition}

\newtheorem{definition}[table]{Definition}

\newtheorem{notation}[table]{Notation}

\newtheorem{convention}[table]{Convention}

\newtheorem{observation}[table]{Observation}

\newtheorem{conjecture}[table]{Conjecture}

\theoremstyle{remark}
\newtheorem{fact}[table]{Fact}

\newtheorem{example}[table]{Example}
\newtheorem{counterexample}[table]{Counterexample}
\newtheorem{exercise}[table]{Exercise}

\newtheorem{problem}[table]{Problem}
\newtheorem{histrmks}[table]{Historical remarks}
\newtheorem{remark}[table]{Remark}
\newtheorem{remarks}[table]{Remarks}

\theoremstyle{plain}
\newtheorem*{thm*}{Theorem}
\newtheorem*{theorem*}{Theorem}
\newtheorem*{prop*}{Proposition}
\newtheorem*{proposition*}{Proposition}
\newtheorem*{lemma*}{Lemma}
\newtheorem*{corollary*}{Corollary}
\newtheorem*{cor*}{Corollary}

\theoremstyle{definition}
\newtheorem*{definition*}{Definition}
\newtheorem*{defn*}{Definition}
\newtheorem*{QQ*}{Question}
\newtheorem*{obs*}{Observation}
\newtheorem*{notation*}{Notation}
\newtheorem*{discussion*}{Discussion}

\theoremstyle{remark}
\newtheorem*{rmk*}{Remark}
\newtheorem*{remark*}{Remark}
\newtheorem*{examples*}{Examples}
\newtheorem*{example*}{Example}
\newtheorem*{EG*}{Example}
\newtheorem*{EGs*}{Examples}
\newtheorem*{fact*}{Fact}
\newtheorem*{prob*}{Problem}

\newcommand{\bthm}{\begin{theorem}}
\newcommand{\ethm}{\end{theorem}}
\newcommand{\bprop}{\begin{proposition}}
\newcommand{\eprop}{\end{proposition}}
\newcommand{\bcor}{\begin{corollary}}
\newcommand{\ecor}{\end{corollary}}
\newcommand{\bconj}{\begin{conjecture}}
\newcommand{\econj}{\end{conjecture}}
\newcommand{\blem}{\begin{lemma}}
\newcommand{\elem}{\end{lemma}}
\newcommand{\bclm}{\begin{claim}}
\newcommand{\eclm}{\end{claim}}
\newcommand{\bpf}{\begin{proof}}
\newcommand{\epf}{\end{proof}}
\newcommand{\bdetails}{\begin{details}}
\newcommand{\edetails}{\end{details}}
\newcommand{\bdefi}{\begin{definition}}
\newcommand{\edefi}{\end{definition}}
\newcommand{\bdefn}{\begin{definition}}
\newcommand{\edefn}{\end{definition}}
\newcommand{\bex}{\begin{example}}
\newcommand{\eex}{\end{example}}
\newcommand{\bprob}{\begin{problem}}
\newcommand{\eprob}{\end{problem}}
\newcommand{\bob}{\begin{observation}}
\newcommand{\eob}{\end{observation}}
\newcommand{\bexer}{\begin{exercise}}
\newcommand{\eexer}{\end{exercise}}
\newcommand{\bexers}{\begin{exercises}}
\newcommand{\eexers}{\end{exercises}}
\newcommand{\brmk}{\begin{remark}}
\newcommand{\ermk}{\end{remark}}
\newcommand{\bhist}{\begin{histrmks}}
\newcommand{\ehist}{\end{histrmks}}
\newcommand{\brmks}{\begin{remarks}}
\newcommand{\ermks}{\end{remarks}}
\newcommand{\bverb}{\begin{verbatim}}
\newcommand{\everb}{\end{verbatim}}
\newcommand{\bntn}{\begin{notation}}
\newcommand{\entn}{\end{notation}}
\newcommand{\bfct}{\begin{fact}}
\newcommand{\efct}{\end{fact}}
\newcommand{\bfcts}{\begin{facts}}
\newcommand{\efcts}{\end{facts}}
\newcommand{\benum}{\begin{enumerate}}
\newcommand{\eenum}{\end{enumerate}}
\newcommand{\bitem}{\begin{itemize}}
\newcommand{\eitem}{\end{itemize}}

\newcommand{\presectionskip}{-1.5\baselineskip}
\newcommand{\postsectionskip}{0.3\baselineskip}
\makeatletter
\usepackage{titlesec}
\renewcommand{\section}{\@startsection
  {chapter}{0}{0mm}
  {\presectionskip}
  {\postsectionskip}
  {\sffamily\huge}}
\renewcommand{\section}{\@startsection
  {section}{1}{0mm}
  {\presectionskip}
  {\postsectionskip}
  {\sffamily\LARGE}}
\renewcommand{\subsection}{\@startsection
  {subsection}{2}{0mm}
  {\presectionskip}
  {\postsectionskip}
  {\sffamily\Large}}
\renewcommand{\subsubsection}{\@startsection
  {subsubsection}{3}{0mm}
  {\presectionskip}
  {\postsectionskip}
  {\sffamily\normalsize}}
\renewcommand{\@seccntformat}[1]{\csname the#1\endcsname.\quad}
\newcommand\HUGE{\@setfontsize\Huge{30}{47}} 
\makeatother
\usepackage{titling}
  \titleformat{\chapter}[display]
  {\sffamily\Large}
  {Chapter {\HUGE\normalfont\thechapter}}    
  {1em}
  {\huge}
  \titleformat{\section}
  {\sffamily\Large}
  {\normalfont{\@setfontsize\Large{18}{47}\thesection.}}    
  {1em}
  {}
  \titleformat{\subsection}
  {\sffamily\large}
  {\normalfont{\Large\thesubsection}.}    
  {1em}
  {}
  \titleformat{\subsubsection}
  {\sffamily\normalsize}
  {\normalfont{\large\thesubsubsection}.}    
  {1em}
  {}

\renewcommand		{\SS}				{\textsection}

\newcommand		{\quation}[1]			{\begin{equation} #1 \end{equation}}

\newcommand		{\eqn}[1]			{\begin{align*} #1 \end{align*}}

\def			\SPSB#1#2			{\rlap{\textsuperscript{#1}}\textsubscript{#2}}
\makeatletter
\def			\smallunderbrace#1		{\mathop{\vtop{\m@th\ialign{##\crcr
							   $\hfil\displaystyle{#1}\hfil$\crcr
							   \noalign{\kern3\p@\nointerlineskip}%
							   \tiny\upbracefill\crcr\noalign{\kern3\p@}}}}\limits}
\makeatother

\makeatletter
\newcommand{\subalign}[1]{%
  \vcenter{%
    \Let@ \restore@math@cr \default@tag
    \baselineskip\fontdimen10 \scriptfont\tw@
    \advance\baselineskip\fontdimen12 \scriptfont\tw@
    \lineskip\thr@@\fontdimen8 \scriptfont\thr@@
    \lineskiplimit\lineskip
    \ialign{\hfil$\m@th\scriptstyle##$&$\m@th\scriptstyle{}##$\crcr
      #1\crcr
    }%
  }
}
\makeatother

\makeatletter
\newcommand		{\oset}[3][0ex]			{%
								\raisebox{.175ex}{$%
								  \mathrel{\mathop{#3}\limits^{
								    \vbox to#1{\kern-2\ex@
								    \hbox{$\scriptstyle#2$}\vss}}}
								    $}%
							    }
\makeatother

\usepackage{faktor,xfrac}

\newcommand		{\dsp}		{\displaystyle}
\newcommand		{\nd}		{\noindent}

\newcommand		{\bs}		{\bigskip}

\newcommand		{\mn}		{\mspace{-2mu}}
\newcommand		{\mnn}		{\mspace{-1mu}}

\newcommand		{\ol}			{\overline}
\newcommand		{\os}			{\overset}
\newcommand		{\us}			{\underset}

\newcommand		{\ul}			{\underline}
\newcommand		{\wh}			{\widehat}

\newcommand		{\wt}			{\widetilde}

\newcommand		{\mr}			{\mathrm}
\newcommand		{\bb}			{\mathbb}

\newcommand		{\f}			{\mathfrak}
\newcommand		{\ms}			{\mathscr}
\newcommand		{\sans}			{\mathsf}

\newcommand		{\g}		{\gamma}
\renewcommand		{\d}		{\delta}
\renewcommand		{\epsilon}	{\varepsilon}
\newcommand		{\e}		{\epsilon}

\newcommand		{\h}		{\eta}

\renewcommand		{\l}		{\lambda}

\newcommand		{\s}		{\sigma}
\newcommand		{\es}		{\varsigma}
\newcommand		{\vp}		{\varphi}
\newcommand		{\w}		{\omega}

\newcommand		{\W}		{\Omega}

\newcommand		{\D}		{\Delta}

\DeclareSymbolFont{cmletters}{OT1}{cmr}{m}{n}
\DeclareMathSymbol{\Ups}{\mathalpha}{cmletters}{"7}
\renewcommand		{\Upsilon}	{\Ups}

\newcommand		{\ceq}		{\coloneqq}
\newcommand		{\eqc}		{\eqqcolon}

\let\union\cup%
\renewcommand		{\cup}		{\mspace{-1mu}\smile\mspace{-1mu}}

\newcommand		{\less}		{\setminus}
\newcommand		{\sub}		{\subseteq}

\DeclareRobustCommand	{\lq}		{\text{\reflectbox{$/$}}}		

\makeatletter
\DeclarePairedDelimiterX
			{\pmodx}[1]	{(}{)}{{\operator@font mod}\mkern6mu#1}
					\renewcommand{\pmod}{%
					  \allowbreak
					  \if@display\mkern18mu\else\mkern8mu\fi
						  \pmodx
					}
\makeatother

\renewcommand		{\th}		{^{\mathrm{th}}}

\renewcommand		{\:}		{\colon}
\renewcommand		{\-}		{^{-1}}

\renewcommand		{\o}		{\circ}

\newcommand		{\adj}		{\dashv}

\ExplSyntaxOn
\NewDocumentEnvironment{adjunctions}{O{}}
{
	\cs_set_eq:cN {@arraycr} \farin_arraycr:
	\keys_set:nn { farin/adjunction } { #1 }
	\begin{array}
		{
			@{ \hspace { \dim_eval:n { \l_farin_left_shift_dim + \l_farin_padding_dim } } }
			r
			@{ {\farin_strut:} \l_farin_symbol_tl {} }
			l
			@{ \hspace { \dim_eval:n { \l_farin_right_shift_dim + \l_farin_padding_dim } } }
		}
	}
	{
	\end{array}
}
\keys_define:nn { farin/adjunction }
{
	leftshift       .dim_set:N = \l_farin_left_shift_dim,
	leftshift       .initial:n = 0pt,
	rightshift      .dim_set:N = \l_farin_right_shift_dim,
	rightshift      .initial:n = 0pt,
	padding         .dim_set:N = \l_farin_padding_dim,
	padding         .initial:n = 6pt,
	symbol          .tl_set:N  = \l_farin_symbol_tl,
	symbol          .initial:n = \longrightarrow,
	verticalspacing .dim_set:N  = \l_farin_vertspac_dim,
	verticalspacing .initial:n = {3pt},
}
\cs_new_protected:Npn \farin_strut:
{
	\vrule height \dim_eval:n { \ht\strutbox + 1.2\l_farin_vertspac_dim }
	depth  \dim_eval:n { \dp\strutbox + \l_farin_vertspac_dim }
	width 0pt
}
\makeatletter
\exp_args:NNo \cs_new:Npn \farin_arraycr:
{
	\@arraycr\hline
}
\makeatother
\ExplSyntaxOff

\newcommand		{\dual}		{\mn^\vee}
\renewcommand		{\.}		{\cdot}

\newcommand		{\x}		{\times}
\newcommand		{\xu}[3]	{\smash{{#2}\us{#1}\times{#3}}}

\newcommand		{\oplushigher}	{\mathbin{\raisebox{.85pt}{$\displaystyle\oplus$}}}

\DeclareMathOperator*	{\otimesvariable}{%
			\mathchoice {\raisebox{.85pt}{$\displaystyle\otimes$}}
						{\raisebox{.85pt}{$\otimes$}}
						{\raisebox{0.7pt}{$\scriptstyle\otimes$}}
						{\raisebox{0.2pt}{$\scriptscriptstyle\otimes$}}
						}
\newcommand		{\tensor}	{\otimesvariable}

\newcommand		{\direct}	{\oplushigher}
\newcommand		{\ox}		{\tensor}
\newcommand		{\+}		{\direct}
\newcommand		{\Tensor}	{\bigotimes}
\newcommand		{\Direct}	{\bigoplus}

\newcommand		{\BBB}		{\mathbf{B}}

\newcommand		{\susp}		{\Sigma}

\newcommand		{\bul}		{\bullet}

\DeclareMathOperator	{\gr}		{gr }

\DeclareMathOperator	{\diag}		{diag}

\DeclareMathOperator	{\rk}		{rk }

\DeclareMathOperator	{\coker}	{coker }

\DeclareMathOperator	{\Tor}		{Tor}

\DeclareMathOperator	{\chara}	{char }

\newcommand		{\cupone}	{\mathbin{{\cup}_1}}
\newcommand		{\cuptwo}	{\mathbin{{\cup}_2}}

\DeclareMathOperator	{\Map}		{Map}

\DeclareMathOperator	{\Hom}		{Hom}

\DeclareMathOperator	{\End}		{End }

\newcommand		{\Orth}		{\mr{O}}
\newcommand		{\SO}		{\mr{SO}}

\newcommand		{\U}		{\mr{U}}
\newcommand		{\SU}		{\mr{SU}}

\newcommand		{\Sp}		{\mr{Sp}}

\makeatletter
\newbox\xrat@below
\newbox\xrat@above
\newcommand		{\xrightarrowtail}[2][]	{%
						  \setbox\xrat@below=\hbox{\ensuremath{\scriptstyle #1}}%
						  \setbox\xrat@above=\hbox{\ensuremath{\scriptstyle #2}}%
				  \pgfmathsetlengthmacro{\xrat@len}{max(\wd\xrat@below,\wd\xrat@above)+.6em}%
  						\mathrel{\tikz [>->,baseline=-.55ex]
              					   \draw (0,0) -- node[below=-2pt] {\box\xrat@below}
                            					    node[above=-2pt] {\box\xrat@above}
                    						   (\xrat@len,0) ;}
						}
\newbox\xrat@below
\newbox\xrat@above
\renewcommand		{\xtwoheadrightarrow}[2][]{%
						  \setbox\xrat@below=\hbox{\ensuremath{\scriptstyle #1}}%
						  \setbox\xrat@above=\hbox{\ensuremath{\scriptstyle #2}}%
				  \pgfmathsetlengthmacro{\xrat@len}{max(\wd\xrat@below,\wd\xrat@above)+.6em}%
						 \mathrel{\tikz [->>,baseline=-.55ex]
					                 \draw (0,0) -- node[below=-2pt] {\box\xrat@below}
					                                node[above=-2pt] {\box\xrat@above}
						                       (\xrat@len,0) ;}
		       				}
\makeatother
\newcommand		{\xmono}	{\xrightarrowtail}
\newcommand		{\mono}		{\xmono{\phantom{\ \, }}}

\newcommand		{\xepi}		{\xtwoheadrightarrow}
\newcommand		{\epi}		{\xepi{\phantom{\ \, }}}
\newcommand		{\longepi}	{\xepi[]{\ \ \ \ }}

\newcommand		{\longto} 	{\longrightarrow}
\newcommand		{\lt}		{\longto}
\newcommand		{\xtoo}		{\xrightarrow} 
\newcommand		{\from}		{\leftarrow}
\newcommand		{\longfrom}	{\longleftarrow}

\newcommand		{\lmt}		{\longmapsto}

\newcommand		{\simto}	{\xrightarrow{\sim}}

\newcommand		{\longsimto}	{\os\sim\longto}
\newcommand		{\isoto}	{\longsimto}
\DeclareRobustCommand	{\longsimfrom}	{\os\sim\longfrom}
\newcommand		{\isofrom}	{\longsimfrom}

\newcommand		{\longbij}	{\longleftrightarrow}

\newcommand		{\To}		{\Rightarrow}
\newcommand		{\inc}		{\hookrightarrow}
\newcommand		{\xinc}		{\xhookrightarrow}
\newcommand		{\longinc}	{\xinc[]{\ \ \ \ }}

\newcommand		{\act}		{\,\raisebox{.25ex}{$\curvearrowright$}\,}

\newcommand		{\vertsim}	{\rotatebox{90}{$\sim$}}

\newcommand		{\hmt}		{\simeq}
\newcommand		{\iso}		{\cong}
\newcommand		{\homeo}	{\approx}

\newcommand		{\F}		{\bb F}

\newcommand		{\Z}		{\bb Z}
\newcommand		{\Q}		{\bb Q}
\newcommand		{\R}		{\bb R}
\newcommand		{\C}		{\bb C}

\newcommand		{\RP}		{\bb R \mr P}
\newcommand		{\CP}		{\bb C \mr P}

\newcommand		{\RPi}		{\RP^\infty}
\newcommand		{\Rpi}		{\RPi}

\newcommand		{\A}		{\bb A}

\DeclareMathOperator	{\id}		{id}

\newcommand		{\ADR}		{A_{\mathrm{dR}}}
\newcommand		{\APL}		{A_{\mathrm{PL}}}

\newcommand		{\CGA}		{\textsc{cga}\xspace}
\newcommand		{\CGAs}		{\textsc{cga}s\xspace}
\newcommand		{\DGA}		{\textsc{dga}\xspace}
\newcommand		{\DG}		{\textsc{dg}\xspace}
\newcommand		{\DGAs}		{\textsc{dga}s\xspace}
\newcommand		{\DGC}		{\textsc{dgc}\xspace}
\newcommand		{\DGCs}		{\textsc{dgc}s\xspace}
\newcommand		{\CDGA}		{\textsc{cdga}\xspace}

\newcommand		{\CDGAs}	{\textsc{cdga}s\xspace}

\newcommand		{\HGA}		{\textsc{hga}\xspace}
\newcommand		{\HGAs}		{\textsc{hga}s\xspace}
\newcommand		{\SHM}		{\textsc{shm}\xspace}
\newcommand		{\SHC}		{\textsc{shc}\xspace}
\newcommand		{\Ai}		{$A_\infty$}

\newcommand{\Sing}{\mathrm{Sing}}

\newcommand		{\kk}		{k}

\newcommand		{\fk}		{{\mathfrak k}}

\renewcommand 		{\H}		{H^*}

\newcommand 		{\eA}		{\e_{\! A}}

\newcommand		{\SSS}		
					{\textsc{sss}\xspace}

\usepackage{tracklang}
\usepackage{multirow,accents}
\usepackage{array}\newcolumntype{R}{>{$}l<{$}}
\usepackage{scalerel}[2016/12/29]

\usepackage{chngcntr}
\counterwithin{table}{section}
\counterwithin{equation}{section}

\theoremstyle{definition}
\newtheorem*{layout*}{Outline}

\usepackage{silence}
\colorlet{jlabel}{Crimson!75!Black}		
\colorlet{jred}{Crimson!75!Black}	
\colorlet{nonisored}{Crimson!90!Black}	
\colorlet{jhmt}{Dandelion!80!Indigo}	
\colorlet{jcmp}{Black!35!White}
\colorlet{jblue}{RoyalBlue!80!Black}	
\colorlet{jviolet}{Violet!60!Black}
\colorlet{compromise}{rgb:Navy!85!Indigo,3;white,1}	
\colorlet{jQ2}{rgb:Indigo!90!Black,3;red,1}	
\colorlet{jA}{RoyalBlue!40!Indigo}		
\colorlet{jX3}{BlueViolet!100!Black}
\colorlet{jX2}{BlueViolet!100!Black}
\colorlet{jX2p}{Dandelion!50!Magenta}	
\colorlet{jgreen}{SeaGreen!85}
\colorlet{jcyan}{Cyan!70!Black}

\let		\epsilon	\varepsilon
\newcommand		{\B}		{\mathbf{B}}
\renewcommand	{\W}		{\BW}
\newcommand		{\Y}		{\Upsilon}
\newcommand	{\ot}		{^{\otimes \mnn 2}}

\renewcommand	{\2}		{\ot}

\newcommand{\Homd}{D}

\newcommand		{\Hflat}		{H_\flat{}\!\mn{}^*}
\renewcommand	{\susp}			{s}
\newcommand		{\desusp}		{s^{-1}}
\def				\iter#1#2{#1^{[#2]}}
\newcommand		{\muA}		{\mu_{\mn A}}
\renewcommand	{\C}		{C^*}
\newcommand		{\I}		{I^*}
\newcommand		{\Dpra}		{\pi_0 \x 0}
\newcommand		{\Dprb}		{\pi_1 \x 0 = 0 \x \pi_0}
\newcommand		{\Dprc}		{0 \x \pi_1}
\newcommand		{\Miso}		{\Xi}
\newcommand		{\lho}			{\l_H^{(1)}}

\newcommand		{\Algs}		{\sans{DGA}}
\newcommand		{\Coalgs}		{\sans{DGC}}
\newcommand		{\GM}		{\sans{Mod}}
\newcommand		{\Tw}		{\sans{T{\mn}w}}

\newcommand	{\EMSS}		{\textsc{emss}\xspace}
\newcommand		{\SHCA}		{\SHC-algebra\xspace}
\newcommand		{\SHCAs}		{{\SHCA}s\xspace}
\newcommand		{\WHC}		{\textsc{whc}\xspace}
\newcommand		{\WHCA}		{\WHC-algebra\xspace}

\DeclareMathOperator*{\T}{%
	\mathchoice {\raisebox{.85pt}{$\displaystyle\underline{{\otimes}}$}}
	{\mathbin{{\raisebox{.85pt}{$\underline{{\otimes}}$}}}}
	{\mathbin{{\raisebox{.7pt}{$\scriptstyle\underline{{\otimes}}$}}}}
	{\mathbin{{\raisebox{.2pt}{$\scriptscriptstyle\underline{{\otimes}}$}}}}
}

\DeclareFontFamily{U}{wncy}{}
\DeclareFontShape{U}{wncy}{m}{n}{<->wncyr10}{}
\DeclareSymbolFont{mcy}{U}{wncy}{m}{n}
\DeclareMathSymbol{\Sha}{\mathord}{mcy}{"58}

\def				\BW			{\boldsymbol{\Omega}}

\newcommand{\Tr}[1]{	\us{#1}	{\mathrm{Tor}}	}

\newcommand{\BA}{\B A}

\def\kk{k}

\newcommand{\simpsetvar}{X}
\newcommand{\simpset}{\simpsetvar_\bul}

\def\susp{s}
\def\desusp{s^{-1}}

\newcommand		{\twist}	{_}

\renewcommand{\Tr}[1]{	\us{#1}	{\mathrm{Tor}}	}

\DeclareMathVersion{normal2}

\begin{document}

\title
	{
		The 
		cohomology of homogeneous spaces
		in historical context
	}
\author{Jeffrey D. Carlson}

\maketitle

\begin{abstract}
	The real singular cohomology ring of a homogeneous space $G/K$,
	interpreted as the Borel equivariant cohomology $H^*_K(G;\R)$,
	was historically the first computation of equivariant cohomology 
	of any nontrivial connected group action.
	After early approaches using 
	the Cartan model for equivariant cohomology
	with $\R$ coefficients	
	and the Serre spectral sequence,
	post-1962 work computing the groups and rings 
	$H^*(G/K;\kk)$ and $H^*_H(G/K;\kk)$
	with more general coefficient rings
	motivated the development of 
	minimal models in rational homotopy theory,
	the Eilenberg--Moore spectral sequence,
	and $A_\infty$-algebras.
	In this essay, we survey the history of these ideas and the
	associated results.
\end{abstract}


One of the most classical algebraic 
invariants of a continuous group action $G \act X$ is the 
Borel equivariant cohomology $H^*_G(X;\kk)$,
defined as the singular cohomology of the homotopy orbit space 
$\defm{X_G} = (EG \x X)/G$
and studied systematically from 1960.
For $\kk$ the real field $\R$, 
equivariant cohomology already appears in
1950 in Henri Cartan's work computing the real cohomology ring
of a homogeneous space $G/K$ for $G$ a compact, connected Lie group 
and $K$ a closed, connected subgroup.
The determination of the cohomology of a homogeneous space 
is thus the ur-example of a computation of a ring-valued invariant of a
nontrivial connected group action.
It was at the same time a motivating example for minimal models in 
rational homotopy theory.

Generalization of Cartan's result to more general coefficient rings
has been similarly fruitful.
Such work 
directly motivated
differential homological algebra
and
the Eilenberg--Moore spectral sequence,
and led to substantial development 
in the field of \Ai- and other up-to-higher-homotopy algebraic structures.
This program, which has lasted seventy-five years,
is perhaps only now nearing its conclusion.
In this article, 
we will relate the history of these generalizations
and the vistas in topology, homological algebra,
and homotopy theory that they opened.\footnote{\
	There is much to discuss.
	This article has been pruned to the satisfy the length restrictions
	of a proceedings volume
	and eventually a lengthier update should appear.
}

In principle, our exposition will assume
algebraic topology and homological algebra only up to the basic properties 
of the Serre spectral sequence (\textcolor{RoyalBlue}{\SSS}),
the spectral sequence of a filtered cochain complex,
the universal principal $G$-bundle $\defm\es = \defm{\es_G}\: EG \to BG$,
and the 
$\defm {\Tor_A(X,Y)} = 
\Direct_{p \geq 0} \Tor^{-p}_A(X,Y)$ 
associated to a pair $X \from A \to Y$ of maps of commutative graded algebras%
\footnote{\ 
	equipped with a standard grading to be reviewed later
}. Comfort with coalgebras and specifically the bar construction
will help with later parts of the story but is not strictly necessary 
for those willing to take some details on faith.

\begin{convention}
We always write $\defm G$ for a compact, connected Lie group,
and $\defm H$ and $\defm K$ for closed, connected subgroups.
All algebraic objects are cochain complexes 
(differential $\defm d$ of degree $1$\footnote{\ 
	\NB: This is not the case in the primary literature, 
	even in cases motivated primarily by singular cohomology.
})
over a commutative base ring $\defm\kk$ with unity,
usually suppressed in the notation $\defm{\ox} = {\ox_\kk}$
for the tensor product
and $\defm {\H X} = \H(X;\kk)$
for the singular cohomology ring.
Graded modules are cochain complexes with $d = 0$.
A \defd{quasi-isomorphism}
is a cochain map inducing an isomorphism in cohomology.
Quasi-isomorphisms induce an equivalence relation on differential
graded algebras (allowing zig-zags of maps in alternating directions), 
and if $A'$ and $A$ lie in the same equivalence class,
we say $A'$ is a \defd{model} of~$A$ (and vice versa). 
A \defd{model} of a space $X$ is a differential graded algebra 
(\textcolor{RoyalBlue}{\DGA}) modeling $A$ the cochain
algebra $\defm{\C}(X) = \C(X;\kk)$, so that $\H(A) \iso \H(X)$.
A smooth manifold $X$ is also modeled by its 
\emph{de Rham algebra} $\defm{\ADR(X)}$ for $\kk = \R$.
\end{convention}

The most optimistic guess for a generalization of Cartan's result
forms a template which we will adapt as we encounter the actual results:

\begin{theorem}[one-sided template]\label{thm:Cartan}
Whenever $\kk$ is chosen such that
$\H(BG)$ and $\H(BK)$ are both polynomial rings,
there is an isomorphism 
of graded $\kk$-algebras
\[
\H(G/K;\kk) 
	\isoto 
\Tor_{\H( BG)}(\kk, \H BK)
\mathrlap.
\]
\end{theorem}

\nd
We call this ``one-sided'' because it corresponds to the right action of $K$
on $G$. 
Given
another closed, connected subgroup $H$, 
there is also a ``two-sided'' action of $H \x K$ on $G$
by $(h,x)\.g = hgx\-$,
leading to a more general guess:

\begin{theorem}[two-sided template]\label{thm:Kapovitch}
Whenever $\kk$ is chosen such that
$\H(BG)$, $\H(BH)$, and $\H(BK)$ are all polynomial rings,
there is an isomorphism 
of graded $\kk$-algebras
\[
\H_H(G/K;\kk) 
	\isoto 
\Tor_{\H (BG)}(\H BH, \H BK)
\mathrlap.
\]
\end{theorem}

\nd
The actual theorems specialize $\kk$ in some way,
show only an additive isomorphism, 
come with some restriction on the cochain algebras,
or weaken the condition on $\H(BK)$,
but we will see Cartan's progenitor
is exactly \Cref{thm:Cartan} for $\kk = \R$.

Most subsequent sections will sketch a proof of a
variant of \Cref{thm:Cartan} or \ref{thm:Kapovitch}.
In broad overview, these were first established for $\kk$
a field of characteristic~$0$ by Cartan and Borel
using the Serre spectral sequence and commutative models.
In the 1970s they were extended, additively,
to more general~$\kk$,
with complications in characteristic $2$
using the Eilenberg--Moore spectral sequence and \Ai-algebraic techniques.
Multiplicative results for $\kk$ not containing~$\Q$ have come only since 2019,
and require $2$ to be a unit.
The most general multiplicative result is the author's 2021 \cref{thm:CF}
joint with Matthias Franz,
while the most general additive result is still Hans J. Munkholm's
1974 \cref{thm:Munkholm}.

\begin{convention}\label{def:convention}
The degree of a homogeneous element $x$ is written $\defm{|x|}$.
Maps $f\: C \to A$
of graded $\kk$-modules 
are all $\kk$-linear,
shifting grading by a fixed degree~$\defm{|f|}$.
We write $f \in \defm{\GM_{|f|}}(C,A)$.
We regard the direct sum $\GM(C,A) = \Direct_{n \in \Z} \GM_n(C,A)$ 
as a cochain complex under the differential
$\defm{Df} \ceq d\mn_A f - (-1)^{|f|}f\mnn d_C$.
A cochain map inducing an isomorphism in cohomology is a
\defd{quasi-isomorphism}.

A \textcolor{RoyalBlue}{\DGA}
is a differential graded $\kk$-algebra $A$
equipped with an augmentation $\e\: A \lt \kk$,
with kernel the augmentation ideal $\defm{\ol A} = \ker \e$.
For a cochain algebra $\C(X)$,
an augmentation is induced by restriction to a basepoint in~$X$.
Write $\defm\Algs$ for the category of \DGAs 
and augmentation-preserving \DGA maps.
The Koszul sign convention 
$(f \ox g)(a \ox b) = (-1)^{|g||a|} fa \ox gb$
is always in effect.
\defd{Commutativity} means graded-commutativity 
$ab = (-1)^{|a||b|}ba$.
A commutative \DGA is a 
\textcolor{RoyalBlue}{\CDGA}.
Commutative graded algebras
(\textcolor{RoyalBlue}{\CGA}s)
are \CDGAs with $d = 0$.
A \DGA (or cochain complex) \defd{computes} a graded $\kk$-algebra (or module)~$H$ if its cohomology is $H$.

A \DGA $A$ comes with an augmentation ideal $\defm{\ol A} = \ker \e$.
Its \defd{ideal of decomposable elements} 
is $\ol A \ol A = \{\sum a_j b_j : a_j,b_j \in \ol A\}$
and its \defd{module of indecomposables}
is $\defm{Q(A)} \ceq \ol A / \ol A\ol A$.
Given an oddly graded vector space $V$ of finite type, there is a natural
pairing between the exterior algebra $\defm{\Lambda[V]}$
and the exterior algebra $\Lambda[V^*]$ on the dual $\defm{V^*} \ceq \Hom_\kk(V,\kk)$.
The \defd{primitives} of $\Lambda[V^*]$
are those elements vanishing on the decomposables $\ol{\Lambda[V]}\.\ol{\Lambda[V]}$
of $\Lambda[V]$.
\end{convention}

\section{Cartan}\label{sec:Cartan}

Cartan's original result from his 1950 announcement~\cite{cartan1950transgression}
interprets $G_K$
as the total space of the Borel fibration $G \to G_K \to BK$
and computes $\H(G/K) = \H(G_K)$ over $\kk = \R$.

A key tool is an acyclic \CDGA $\defm{W(\fk)}$ 
equipped with 
certain operations of the Lie algebra~$\fk$ of $K$,
due to Weil (unpublished)
and now called the \emph{Weil algebra},
which serves as a universal model for the data associated with a connection
on a principal $K$-bundle.
By construction,
a connection on a principal $K$-bundle $\defm\pi\: E \to B$
corresponds to a unique \defd{characteristic}
\DGA \defd{map} $ W(K) \lt \ADR(E)$
preserving the operations of $\fk$.
The Weil algebra can also be understood as a model for the total space of 
the universal principal $K$-bundle $\varsigma\: EK \to BK$.
Although $BK$ was only defined in full generality by Milnor later (1956),
it was understood from known cases that $\H(BK)$
should be isomorphic to the invariant subring $S[\f k\dual]^K$ 
of the symmetric algebra $S[\f k\dual]$ 
on the dual $\f k\dual$ 
under the coadjoint $K$-action, graded with $|\f k| = 2$.
The \CGA underlying the Weil algebra is the tensor product
of $S[\fk\dual]$
and the exterior algebra $\Lambda[\f k\dual]$,
graded with $|\fk\dual| = 1$,
and the inclusion $ S[\f k\dual]^K \lt W(K)$ 
can be seen as a model of the universal bundle projection $\varsigma$.
The characteristic map associated to a connection 
restricts to a map
$S[\f k\dual]^K \to \pi^*\ADR(B) \simto \ADR(B)$,
and the induced map in cohomology is
the Chern--Weil homomorphism $\defm{\chi^*}\: S[\f k\dual]^K \lt \H(B)$.

Given a principal $K$-bundle $E \to B$, 
Cartan views $W(K) \otimes \ADR(E)$ 
as a model for forms on $EK \x E \hmt E$.
The projection 
to $E_K = (EK \x E)/K \hmt E/K = B$
induces a pullback 
identifying forms on $E_K$ with
a subalgebra $\defm {\wh C}$ of $W(K) \otimes \ADR(E)$
called the \emph{Weil model}. 
The projection 
$W(K) = S[\f k\dual] \otimes \Lambda[\f k\dual] \otimes \ADR(E)
\lt 
S[\f k] \otimes \ADR(E)$ 
induces an algebra isomorphism between $\wh C$ and
$\defm{C} = \big(S[\f k\dual] \otimes \ADR(E)\big){}^K$.
Transferring the differential along this isomorphism
makes $C$ a \DGA called the \emph{Cartan model}
computing $\H(E_K) = \H_K(E)$.

Results of Hirsch allow one to define 
a subcomplex $C'$ of $C$,
isomorphic to 
$S[\f k]^K \otimes \H(E)$,
such that the inclusion $C' \inc C$ is a quasi-isomorphism.
In the special case  $\pi\: E \to B$ is 
the quotient projection $G \to G/K$,
one can take the ring $\ADR(G)^G$
(isomorphic to $\Lambda[\f g\dual]^G \iso \H(G)$)
as a set of representatives of $\H(G)$,
so that Cartan obtains a \DGA structure 
on $S[\f k\dual]^K \otimes \Lambda[\f g\dual]^G$
computing $\H(G/K)$ as a ring.

\bthm\label{thm:Cartan-detail}
One has the commutative diagram of \CGAs  
\[
\xymatrix@R=0em@C=1.5em{
&\H\big(S[\f k\dual]^K \otimes \Lambda[\f g\dual]^G\big) \mathrlap{ { } = \H(C')} 
                        \ar[dd]^\vertsim \ar@{->>}[dr]&\\
\                        S[\f k\dual]^K\ \ar@{^{(}->}[ur]\ar[dr]_{\chi^*}& & 
\Lambda[\f g\dual]^G \mathrlap{{ } = \H(G).}\\
&\H(G/K)\ar[ur]_{\pi^*}&
}
\]
\ethm

Cartan explores the sequence
$S[\f g\dual]^G \to S[\f k\dual]^K \to \H(G/K) \to \H(G) \to \H(K)$
and as an example computes the Poincar\'e polynomials
of the real oriented Grassmannians
$\wt G_\ell(\R^{\ell+m}) = \SO(\ell+m)/\big(\SO(\ell) \x \SO(m)\big)$.
He also notes that his result implies
$\H(G/K) \iso S[\f k\dual]^K \ox_{S[\f g\dual]^G} \R$
as rings
when $K$ is of full rank in $G$.\footnote{\ 
    This is also a result of Leray
    which had already been published in the case $G$ is finitely covered
    by a product of classical groups~\cite{lerayCR1949e,leray1950c}.
}

\bex
Taking $G = \SO(2n+1)$ and $K = \SO(2) \x \SO(2n-1)$ for $n > 1$
and using restriction relations between Pontrjagin and Euler classes,
one finds $\H \wt G_2(\R^{2n+1}) \iso \R[e]/(e^{2n})$
for $e$ the image of the universal Euler class in $H^2 B\SO(2)$.
\eex

The differential of $C'$
in 
\Cref{thm:Cartan-detail}
is the unique derivation vanishing on~$S[\fk\dual]^K$
and extending a certain linear map 
on a space of exterior generators $\defm{T_G}$ of~$\Lambda[\f g\dual]^G$.
Namely, 
$d|_{T_G}: T_G \to S[\fk\dual]^K$ is the composite of the restriction 
$S[\f g\dual]^G \to S[\f k\dual]^K$
and a map $\wt\tau\: T_G \to S[\f g\dual]^G$
called a \emph{(choice of) transgression}.

In a cohomological first-quadrant spectral sequence $E_*^{*,*}$,
a class of $\smash{E_2^{0,p-1}}$
is said to \defd{transgress} if it survives to $\smash{E_{p}^{0,p-1}}$.
Then $d_{p}$ is defined on the class $[z]_{p} \in E_{p}^{0,p-1}$,
determining an image $\defm\tau[z] = d_{p}[z]_{p} \in E_{p}^{p,0}$
called its \defd{transgression}.
Thus $\tau$ is a degree-$1$ linear map
from a submodule of $\smash{E_2^{0,*}}$ to a quotient of $\smash{E_2^{*,0}}$.
It is often noncanonically lifted to a map $\defm{\wt\tau}$
to $\smash{E_2^{*,0}}$ called a \defd{choice of transgression}.
Cartan's $\wt\tau$ is that of an algebraic 
spectral sequence modeled on the \SSS of $K \to EK \to BK$,
converging from $\H(BK) \otimes \H(K) \iso 
S[\fk\dual]^K \otimes \Lambda[\fk\dual]^K$
to $\H\big(W(K)^K\big) = \H(EK) = \R$.\footnote{\ 
	He actually avoids mentioning spectral sequences as follows.
	Because the Weil algebra $W(K)$ is acyclic and 
	$\iota\: S[\fk\dual]^K \inc S[\fk\dual] \otimes \Lambda[\fk\dual] \iso W(K)$
	is a \DGA map,
	each $x \in S[\fk\dual]^K$ must be $d_{W(K)}y$ for some $y \in W(K)^K$,
	which the projection to $\Lambda[\fk\dual]$
	takes to some $z \in \Lambda[\fk\dual]^K$.
	The resulting assignment $\defm\s\:
	S[\fk\dual]^K \lt \Lambda[\fk\dual]^K$
	taking
	$x$ to $z$ is easily seen to be well-defined,
	and its image is the space ${T_K}$ of transgressive elements.
	Cartan's choice of transgression
	is any linear section $\wt\tau$ of $\s$.
	}
\begin{theorem}[Cartan--Chevalley]\label{thm:Cartan-transgression}
The subspace $\defm{T_K}$ of transgressive elements
in $\Lambda[\fk\dual]^K$
is the space of primitives (\Cref{def:convention}).
A basis $(z_j)$ of $T_K$ forms an irredundant set of exterior generators
of $\Lambda[\fk\dual]^K$.
The codomain of the transgression is the space
of indecomposable elements,
$S[\fk\dual]^K/ S^{\geq 1}[\fk\dual]^KS^{\geq 1}[\fk\dual]^K$.
Under any choice of transgression $\wt\tau$, 
the images $\wt\tau z_j$ form an irredundant set 
of polynomial generators for $S[\fk\dual]^K$.\footnote{\
	Cartan does not include a proof,
	and notes that this work is inspired in part from Koszul's thesis,
	which defines Lie algebra cohomology
	and studies the transgression in 
	a spectral sequence analogous to the \SSS of $K \to G \to G/K$,
	and answers a May 1949 conjecture of Weil.}
\end{theorem}
Thus a choice of transgression $\wt\tau$ induces a linear bijection
between exterior generators of 
$\Lambda[\fk\dual]^K \iso \H(K)$
and polynomial generators of
$S[\fk\dual]^K \iso \H(BK)$.

\section{Borel}\label{sec:Borel}
The proof of \Cref{thm:Cartan-transgression}
relies heavily on the structure of real Lie algebras.
Borel's 1952 thesis, among other things, 
generalizes this transgression theorem as a result about 
spectral sequences over other base fields $\kk$.

\begin{theorem}\label{thm:Borel-trans}
	Let $\kk$ be a field and
	$E_*^{*,*}$ a first-quadrant cohomological spectral sequence
	of bigraded algebras
	such that $E_2^{*,*} \iso E_2^{*,0} \otimes E_2^{0,*}$
	as a bigraded algebra
		and $E_\infty^{*,*} = E_\infty^{0,0} = \kk$.
	\bitem
	\item
		If $E_2^{0,*}$ is an exterior algebra on odd-degree elements,
		then there exist 
		homogeneous transgressive $z_j$
		such that $E_2^{0,*} = \Lambda[z_j]$.
	\item
		If $\chara \kk = 2$,
		suppose there exist
		homogeneous transgressive $z_1,\ldots,z_n \in E_2^{0,*}$ 
		(of any degree)
		such that the  monomials $z_{j_1}\cdots z_{j_\ell}$
		($\mn j_1 < \cdots < j_\ell$, $\ell \leq n$)
		form a basis of $E_2^{0,*}$.
	\eitem
	In either case, for any choice $\wt\tau$ of transgression,
	$E_2^{*,0} = \kk[\wt\tau z_j]$.\footnote{\ 
		This can be strengthened by requiring only 
		$\smash{\Direct_{p+q = 1}^n E^{p,q}_\infty} = 0$
		and concluding only $E_2^{\leq n,0}$ is polynomial,
		where $\max_j |z_j| \leq \frac n 2 - 1$.
	} 
\end{theorem}

%

\textbf{In the rest of this section 
we consider a bundle
$\smash{F \os {\defm i}\to E \os{\defm\varpi}\to B}$
with trivial $\pi_1(B)$-action on $\H(F)$.}
In the \SSS of this bundle,
the transgression goes from $\H(F)$ to $\H(B)$.\footnote{\ 
    Unpacking the definition,
    if a cocycle $z \in C^{p-1}(F)$ represents
    a transgressive class $[z] \in H^{p-1}(F) = \smash{E_2^{0,p-1}}$,
    then there exist a cochain $y \in C^{p-1}(E)$ with $i^*y = z$
    and a cocycle $x \in C^p(B) = \smash{E_2^{p,0}}$ 
    such that $\pi^*x = \d y$,
    and then $d_p[z]_p = [x]_p$.
}
Applied to $K \to EK \os{\es}\to BK$,
\Cref{thm:Borel-trans} 
says that if~$\H(K)$ is exterior over a field $\kk$ of characteristic $\neq 2$,
then $\H(BK)$ is polynomial on a basis of transgressions.\footnote{\ 
	The degree-truncated version of Borel's transgression
	theorem is relevant because 
	Borel and Cartan used only finite-dimensional truncations of $BG$,
	which would later be defined for general topological groups $G$ by Milnor.
} 
In retrospect, this proves \Cref{thm:Cartan} for $K = 1$ and $G/K = G$;
in dealing with $\H(G)$,
Borel is also the first to characterize finitely-generated 
commutative Hopf algebras over $\F_p$.
He also shows that in the \SSS of $\varsigma$,
if $\kk$ is such that $\H(K)$ is a free module on monomials
in transgressive generators as in \Cref{thm:Borel-trans},
the map $\H(K) \to \H(K \x K) \simto \H(K) \ox \H(K)$ induced by 
the group multiplication $K \x K \to K$ takes $z$ to $1 \ox z + z \ox 1$
precisely for these transgressive generators.
These elements are also called \defd{primitive},
and this agrees with the previous notion if $\H(K)$ is exterior.

To compute $\H(G/K)$,
we want a functorial $\R$-\CDGA $\defm A(-)$
computing cohomology of spaces.
	The de Rham algebra $\ADR(-)$ does this only for manifolds.\footnote{\
Borel extends this to compact separable metric spaces $X$
of finite dimension,
using the Menger--N\"obeling theorem
asserting a homeomorphic embedding $X \mono \R^{1 + 2 \dim X}$,
	 restricting the sheaf $U \mapsto \ADR(U)$ on $\R^{1 + 2 \dim X}$ to $X$
	 and taking global sections.
	 }
Sullivan later introduced the $\Q$-algebra $\APL$ of polynomial
differential forms,
and $\kk \ox_\Q \APL(-)$
works for all fields $\kk$ of characteristic $0$.
Equipped with such a model, Borel can use the \SSS  
to generalize Cartan's model $C'$ from
\Cref{thm:Cartan}.
\textbf{From now on, suppose additionally
$\H(F)$ is exterior on a set of generators transgressing in the} \SSS\textbf{.}
If $F = G$ and $\varpi$ is a principal $G$-bundle,
\Cref{thm:Borel-trans}
shows this happens if $\chara\kk \neq 2$
or $\chara\kk = 2$ and the generators $z_j$ transgress in the \SSS
of the universal bundle:
indeed, the classifying map $\chi$ 
induces a map of {\SSS}s from that of $\varsigma$ to that of $\varpi$,
so that in the latter,
each $z_j$ transgresses to $\tau_\varpi (z_j) = \chi^*\tau_{\es} (z_j) \in \H(B)$.

Borel endows the graded algebra $\defm L = A(B) \ox \H(F)$
with the unique derivation extending that on $A(B) \iso A(B) \ox \R$
and taking each $1 \ox z_j$ 
to a cocycle in $A(B) \ox \R$ representing 
$\tau z_j \ox 1 \in \H(B) \ox \R$.
For each~$z_j$, one can find an $A(E)$-cochain~$y_j$
restricting along $i$ to representative of $z_j$ in $A(F)$,
and such a choice
induces an algebra map $j\: \H(F) \to A(E)$ by the commutativity of $A(E)$.
There is thus a \DGA map 
$L \lt A(E)$
taking
$b \ox y \lmt \varpi^*(b)\.j(y)$.

\bthm[{\cite[Thm.~24.1$'$]{borelthesis}}]\label{thm:Borel-Chevalley}
Under the \DGA structure on $L = A(B) \ox \H(F)$
induced by the transgression as in the previous paragraph,
the map $L \to A(E)$ induces a $\kk$-\CGA isomorphism
$\H\big(A(B) \ox \H(F)\big) \isoto \H(E).$
\ethm
\bpf
Replacing $B$ by a weakly equivalent CW complex if necessary,
$A(E)$ inherits the Serre filtration
$F_p A(E) = \ker\big(A(E) \to A(\varpi^*B_{p-1})\big)$
and $A(B) \ox \H(F)$ the filtration 
by $\ker\!\big(A(B) \to A(B_{p-1})\big) \ox \H(F)$.\footnote{\ 
	This is a simplification;
	Borel's actual filtration is essentially by degree of forms.
}
The map $L \to A(E)$
respects these filtrations,
inducing a map of spectral sequences
which reduces on $E_2$ pages to $\id_{\H(B) \ox \H(F)}$
and hence is a quasi-isomorphism.
\epf

\textbf{Suppose additionally from now on
that there exists a sub-}\DGA \textbf{$\defm H$ of $A(B)$
such that the inclusion is a quasi-isomorphism.}
Then one can select representatives of 
$\tau z \ox 1 \in \H(B) \ox \R$ lying in $H \ox \R$
to define a sub-\DGA $\defm{L'} = H \ox \H(F)$ of $L$,
and another filtration spectral sequence argument
shows $L' \inc L$ is a quasi-isomorphism,
so the composite $L' \to A(E)$ is as well.
\begin{theorem}%
	[{\cite[Thm.~25.1]{borelthesis}}]\label{thm:Borel-Cartan}
\!\!\footnote{\ 
	Borel restricts to principal $G$-bundles 
	in his statement,
	using his Thm.~24.1, but using Thm.~24.1$'$
	instead (\Cref{thm:Borel-Chevalley}),
	the restriction is unnecessary.
}
Under the \DGA structure on $L'$
induced by the transgression as in the previous paragraphs,
the map $L \to A(E)$ induces a $\kk$-\CGA isomorphism
$\H\big(\H(B) \ox \H(F)\big) \isoto \H(E).$
\end{theorem}

The hypotheses of \Cref{thm:Borel-Cartan} are always satisfied
for $E = F_K \lt  BK = B$
the Borel fibration associated to the action
of a compact, connected Lie group $K$-action on $F$:
any choice of $A$-cocycles representing generators of the polynomial ring $\H(BK)$
induces a unique \DGA map $\H(BK) \lt A(BK)$,
and we may take for $H$ its image,
yielding an isomorphism
$\H\big(\H(BK) \ox \H(F)\big) \lt \H(F_K)$.
Taking $F = G$ a compact Lie group containing $K$,
with the translation action, 
we get back Cartan's isomorphism
$\H\big(\H(BK) \ox \H(G)\big) \lt \H(G/K)$ from \Cref{thm:Cartan}%
~\cite[Thm.~25.2]{borelthesis}.

Borel is able to obtain results when $\rk K = \rk G$
extending Leray's results over a field of characteristic $p$ or $\Z$:
if $T$ is a maximal torus of $K$ and
the cohomology of $G$, $K$, $G/T$, and $K/T$
are assumed to be $p$-torsion--free for characteristic $p$
or torsion-free for $\kk = \Z$,
then $\H(G/K) \iso \H(BK) \ox_{\H(BG)} \kk$.
He provides several explicit computations,
applying the equal-rank result to the cohomology rings of
$\U(\sum n_j)/\prod \U(n_j)$ and
 $\Sp(\sum n_j)/\prod \Sp(n_j)$ for $\kk = \Z$
 and to $\SO(2n)/\U(n)$ for $\chara\kk \neq 2$.
These computations use Chevalley's restriction isomorphisms
$\H(BG) \to \H(BT)^{W_G}$
for $G$ and $K$,
for $W_G$ the Weyl group,
which holds assuming the torsion conditions are satisfied.
For $\chara\kk =2$, he shows $\SO(2n) \epi \SO(2n)/\U(n)$
induces an injection in cohomology onto the subalgebra
generated by even-degree primitives.

Other spaces treated are $\U(2n)/\Sp(n)$ for $\kk = \Z$
and $\U(n)/\SO(n)$ for $\chara\kk \neq 2$.
In later work he deals with 
the same space and with
$\Orth(\sum n_j)/\prod \Orth(n_j)$ 
for $\chara\kk = 2$,
using instead of a torus $T$
the diagonal elementary abelian $2$-group.
Using the already-studied \SSS of $K \to G \to G/K$ instead,
he also analyzes the Stiefel manifolds $\U(n)/\U(\ell)$ 
and $\Sp(n)/\Sp(\ell)$ for $\ell < n$ and $\kk = \Z$,
and $\SO(n)/\SO(\ell)$ 
in characteristic $2$
and additively for $\kk = \Z$.

\bex
Taking $G = \U(n)$ and $K = \U(1)^n$ diagonal,
one finds 
\[
\H(G/K;\Z) \iso \Z[t_1,\ldots,t_n]/
\big(\sum t_j, 
\smash{\sum_{i < j} t_i t_j}, 
\smash{\sum_{h < i < j} t_h t_i t_j},
\cdots, t_1 \cdots t_n\big)
\mathrlap,
\]
where $|t_j| = 2$.
\eex

\section{Eschenburg}\label{sec:Eschenburg}

There has long been a school of geometers interested 
in compact Riemannian manifolds of positive curvature.
There seem to be relatively few of these,
and many known examples are \defd{biquotients},
orbit spaces of a Lie group $G$
under a free ``two-sided'' action $(u,v)\.g = ugv^{-1}$
of a subgroup $U$ of $G\x G$.
This obviously specializes to a homogeneous space for $U = 1 \x K$. 
The special case where $U = H \x K$ for two closed subgroups $H,K \leq G$
is written
$\defm{H\lq G / K}$,
and this case 
is in fact general,
since writing $\D\: G \to G \x G$ for the diagonal,
$(x,y) \mapsto xy\-\: G \x G \to G$
induces a natural diffeomorphism
$G/U = U \lq (G \x G) /\D G$.

Eschenburg~\cite{eschenburg1992biquotient}
studied the cohomology of biquotients, 
showing many aspects of Borel's and Cartan's analyses generalize.
For convenience, given a right $G$-space $X$ and left $G$-space $Y$,
write $\defm{X \ox_G Y}$ for the orbit space under $(xg,y) \sim (x,gy)$.
Consider the two-sided action of $G \x G$ on $G$ and the restricted
action of~$U$.
Eschenburg notes the following system of quotients of $EG \x EG \x G$:%
\footnote{\ 
The only non-obvious equivalence may be the homeomorphism 
$(e,e') \ox x \mapsto (ex,e')\D G$ on the upper right.
It is well-defined because 
$(eg,e'g') \ox g^{-1}xg' = (e,e') \ox x$
is also sent to $(eg\.g^{-1}xg',e'g')\D G = (exg',e'g')\D G = (ex,e')\D G$.
It is obviously surjective.
It is injective because if $(fy,f')\D G = (ex,e')\D G$,
there is $g \in G$ with $(fyg,f'g) = (ex,e')$,
so $e = fygx^{-1}$,
and $(e,e') \ox x = (fygx^{-1},f'g) \ox x = (fyg,f'g) \ox x^{-1}x1 = 
(fy, f') \ox g^{1-}1g = (f,f') \ox y$.
}
\begin{equation}\label{eq:Eschenburg}
\xymatrix{
G_{U} \hmt (EG \x EG) \ox_{U} G \ar[r]\ar[d]_{\defm\varpi}&
(EG \x EG) \ox_{G \x G} G \homeo (EG \x EG)/\D G \ar[d]_{\defm\d}\\
BU \hmt 
(EG \x EG) \ox_{U} {\ast} \ar[r]_(.425){\defm\chi}&
(EG \x EG) \ox_{G \x G} {\ast} \homeo BG \x BG\mathrlap.
}
\end{equation}
Here the upper-left corner is of interest
because it is homotopy equivalent to $U$
when the action is free.
Since the fiber of both $\varpi$ and $\d$ is $G$,
this is a bundle map, inducing a map of {\SSS}s.
Because $\d$ can be identified with the diagonal map
$BG \lt BG \x BG$ up to homotopy,
$\H(\d)$ can be identified with the cup product
on the polynomial algebra $\H(BG)$,
a surjection with kernel generated by 
the elements
$1 \ox \tau z_j -\tau z_j \ox 1 \in \H(BG) \ox \H(BG) \iso \H(BG \x BG)$
for $z_j$ generators of the exterior algebra $\H(G)$
and $\tau = \tau_\es$ a choice of transgression for $G \to EG \os\es\to BG$.
Since none of the $z_j$ survive the \SSS of $\d$,
it follows each transgresses to $\tau_\d (z_j)= 1 \ox \tau z_j -\tau z_j \ox 1$.
Applying the map to the \SSS of $\varpi$, one finds the following:
\bprop[Eschenburg]
Each $z_j$ transgresses to $\tau_\varpi(z_j) = \chi^* \tau_\d(z_j)$.
\eprop
Eschenberg uses this to study the \SSS of $\varpi$ 
and compute several examples.

%

\section{Kapovitch}\label{sec:Kapovitch}

Eschenburg noted his $G$-bundle map $\varpi \to \d$
implies exterior generators of $\H(G)$ 
transgress in the \SSS of $\varpi$.
Taking $U = H \x K$,
in his thesis work in 2014,
the author noticed that
\Cref{thm:Borel-Cartan},
which is essentially in Borel's thesis,
yields a model
$\H(BH \x BK) \ox \H(G)$ 
for $\H(G_{H \x K}) \iso \H_H(G/K)$.
It turns out this had been known for ten years
for more general reasons.

Borel's quasi-isomorphism $L' \lt A(G_K)$
can be seen
as an early example of rational homotopy theory in action 
(and so, with more squinting, the results of Cartan, Chevalley, Koszul, and Weil
can be as well).
A \emph{Sullivan algebra} is \CDGA
over a field $\kk$ of characteristic $0$
which is free (exterior $\ox$ polynomial) as a \CGA
and whose differential satisfies a certain nilpotence property,
and one constructs and computes with \emph{Sullivan models} 
of spaces,
Sullivan algebras computing their cohomology.
Borel's \cref{thm:Borel-Cartan} 
applies to the universal $G$-bundle $\es\: EG \to BG$
by his transgression \cref{thm:Borel-trans},
defining a model $\H(BG) \ox \H(G)$ for $EG$
with differential the derivation
defined as $0$ on $\H(BG)$ 
and on exterior generators by $z_j \lmt \tau z_j$.
This is a Sullivan model for $EG$.
The bundle $G_K \to BK$ is the pullback of $\es$
under the classifying map $\rho = B(K \inc G)$,
and $\H(\rho)\: \H(BG) \to \H(BK)$ is a 
map of Sullivan models, modeling $\rho$,
and the inclusion $\H(BG) \inc \H(BG) \ox \H(G)$ models $\es$.
A standard result on Sullivan models%
~\cite[\SS15(c)]{FHT} 
says in essence that given a map $X \to B$ to a simply-connected $B$, 
a Serre fibration $E \to B$,
and Sullivan models $M_X \from M_B \inc M_E$,
finitely generated in each degree,
and such that
$M_B \inc M_E$ is a \emph{relative Sullivan algebra},
meaning the map of underlying \CGAs
is $M_B \inc M_B \ox A$ for some other \CGA $A$, 
the tensor product $M_X \ox_{M_B} M_E$ is a Sullivan model for 
the pullback $X \x_B E$.
Applying this, we get back Borel's model for $G_K$:
\quation{\label{eq:Koszul}
\H(BK) \ox_{\H(BG)} \big(\H(BG) \ox \H(G)\big) \iso \H(BK) \ox \H(G) 
= L'\mathrlap.
}

Building on Eschenburg's work,
Kapovitch notes 
\[
\H(\chi)\: \H(BH) \ox \H(BK) \to \H(BG) \ox \H(BG)
\]
is a Sullivan model of $\chi = B(H \x K \inc G \x G)$
and Eschenburg's transgression result 
yields a Sullivan model 
\[
\defm R = \H(BG) \ox \H(G)\ox \H(BG)
\]
of $(EG \x EG)/\D G$
with differential vanishing on $\H(BG) \ox \H(BG)$
and defined by the transgressions 
$\tau_\d(z_j) \lmt 1 \ox \tau z_j - \tau z_j \ox 1$
on exterior generators $z_j$ of $\H(G)$.
Moreover the inclusion $\H(BG) \ox \H(BG) \inc R$ is a model of $\d$.
Then the same result on pullbacks and Sullivan models gives the following.

\bthm[{\cite{kapovitch2002biquotients}}]\label{thm:Kapovitch-detail}
The equivariant cohomology $\H_H(G/K)$ over $\kk = \Q$
is the cohomology of the Sullivan algebra
\[
\big(\H(BH) \ox \H(BK)\big) \ox_{\H(BG)\ox\H(BG)} R 
\iso \H(BH) \ox \H(BK) \ox \H(G) 
\]
with differential the derivation vanishing on $\H(BH) \ox \H(BK)$
and defined on exterior generators of $\H(G)$
by $z_j \lmt 1 \ox \rho_K^*\tau z_j - 1 \ox \rho_H^*\tau z_j$
for $\rho_H = B(H \inc G)$ and $\rho_K = B(K \inc G)$.
\ethm

\bex[{\cite[Prop.~1.9]{carlson2016grassmannian}}]
Let $G = \SU(n)$ and $K \iso \U(1)$  
a \emph{reflected circle},
meaning that for some $g \in G$,
for all $x \in K$,
one has $gxg\- = x\-$.
Then $\H_K(G/K) \iso \Lambda[z_5,\ldots,z_{2n-1}] \otimes \Q[s,t]/(s^2-t^2)$,
where $|s| = |t| = 2$ 
and $\H(G) \iso \Lambda[z_3,z_5,\ldots,z_{2n-1}]$.
\eex

\bex[{\cite{he2016grassmannian}\cite{carlson2016grassmannian}}]
For $G = \SO(\ell + m)$ 
and $K = \SO(\ell) \x \SO(m)$, so that $G/K$ is the oriented Grassmannian,
the ring $\H_K(G/K;\Q)$ has been expressed
in work of the author using the Kapovitch model, 
following a computation by different means by Chen He.
\eex

\section{The Eilenberg--Moore spectral sequence}\label{sec:EMSS}

A critical feature of the approaches 
described so far is the 
use of abstractions of differential forms
to give \CDGA models of spaces
over fields of characteristic $0$.
Functorial \CDGA models 
do not exist in other characteristics,\footnote{\ 
The proof from Borel's
1951 ETH lectures on the Leray spectral sequence is as follows%
~\cite[Thm.~7.1]{borel1951leray}.
	Suppose for a contradiction that 
	$\defm A$ is a $\kk$-{\CDGA}--valued
	contravariant functor on topological spaces,
	for $\kk$ a ring of characteristic $p > 0$,
	such that
	$\H\big(A(-)\big) \iso \H(-)$
	and $i^*\: A(Y) \lt A(X)$ is surjective 
	whenever $i\: X \longinc Y$ is the inclusion of a closed subset.

	We show this impossible 
	for $X = \CP^n$ with $n > p$ and $Y \hmt {*}$ the cone on $X$.
	Note that for any even-degree $y$ in a $\kk$-\CDGA,
	one has $d(y^p) = py^{p-1}dy = 0$.
	Let the cocycle $a \in A^2(X)$
	represent a generator $x$
	of $\H(X) \iso \kk[x]/(x^{n+1})$,
	which since $X$ is closed in $Y$ 
	is $i^* \wt a$ for some $\wt a \in A^2(Y)$.
	Now $\wt a^p \in A^{2p}(Y)$ is a cocycle,
	and since $Y$ is contractible,
	also a coboundary.
	But then as~$i^*$ is a \DGA map,
	$ i^*(\wt a^p)=a^p $ 
	is also a coboundary
	despite representing $x^p \neq 0$ in $H^{2p}(X)$.
}
but one aspect of the computation does generalize.

The model $\H(BG) \ox \H(G)$ of $EG$
from the previous section, regraded by exterior degree,
is a free $H^*(BG)$-module resolution 
of $\kk \iso H^*(BG) / H^{\geq 1}(BG)$,
so Borel's model $\H(BK) \ox \H(G)$ for $G/K$
computes $\Tor_{\H BG}(\kk, \H BK)$.
This finally explains why we have stated Cartan's \cref{thm:Cartan-detail}
as \Cref{thm:Cartan} for $\kk = \R$.
The spaces of interest fit into 
\eqref{fig:Cartan}.

\begin{figure}
\begin{subfigure}{0.3\textwidth}
	\centering
	$
\xymatrix@C=2em{
G_K \ar[r]\ar[d]&		 EG \mathrlap{{ } \hmt {*}}\ar[d] \\
BK \ar[r]_{\rho\vphantom{\rho_H}}		&	BG	
}
$
\subcaption{ }
\label{fig:Cartan}
\end{subfigure}
\hfill
\begin{subfigure}{0.3\textwidth}
	\centering
	$
\xymatrix@C=1em{
(G/K)_H \ar[r]\ar[d]&		 EG/K \mathrlap{{}\hmt BK} \ar[d]^{\rho_K} \\
BH \ar[r]_{\rho_H}		&EG/G \mathrlap{{ }\hmt BG}	
}
$
\subcaption{ }
\label{fig:Kapovitch}
\end{subfigure}
\hfill
\begin{subfigure}{0.3\textwidth}
	\centering
	$
\xymatrix{
\vphantom{X_X} Y \ar[r]\ar[d]_{\defm\varpi}	&	E \ar[d]^{\defm\pi} \\
X \ar[r]_{\chi\vphantom{\rho_H}}		&	B
}
$
\subcaption{ }
\label{fig:genpull}
\end{subfigure}
\caption{Some homotopy pullback squares}
\label{fig:pullbacks}
\end{figure}

We obtain \Cref{thm:Kapovitch} for $\kk = \Q$
from \Cref{thm:Kapovitch-detail} in similar fashion.
The Sullivan model $R$ of $BG$ from the previous section,
viewed as a left $\H(BG)$-module,
yields an (inefficient) $\H(BG)$-module resolution of $\H(BG)$ itself
via the cup product
\[
R_0 = \H(BG) \ox H^0(G) \ox \H(BG) \lt \H(BG)\mathrlap.
\]
Thus $R \ox_{\H(BG)} \H(BK)$
computes $\Tor_{\H(BG)}(\H BG, \H BK) = \H(BK)$,
and so $R \ox_{\H(BG)} \H(BK)$ 
is an $\H(BG)$-module resolution of $\H(BK)$.
Thus 
$\Tor_{\H(BG)}(\H BH, \H BK)$
can be computed as the cohomology of
\[
	\H(BH) \ox_{\H(BG)} R \ox_{\H(BG)} \H(BK)
		\iso 
	\H(BH) \ox H(G) \ox  \H(BK)\mathrlap,
\]
which is the Kapovitch model of $G_{K \x H}$ from \Cref{thm:Kapovitch-detail}.%
\footnote{\
	We promised we would discuss the grading.
	A projective resolution $P_\bul \lt M$ of an $A$-module $M$
	is a sequence of degree-$0$ $A$-module maps $P_p \lt P_{p+1}$ 
	for $p \leq 0$ together with a degree-$0$
	$A$-module map $P_0 \lt M$
	such that the sequence
	$\cdots \to P_{-1} \to P_0 \to M \to 0$
	is exact.
	A pure tensor $x \ox y$ 
	of homogeneous elements $x \in P_p$ and $y \in N$
	inherits a well-defined ``internal'' degree $q = |x| + |y|$
	in addition to the nonpositive resolution degree $p$,
	inducing a bigrading on $P_\bul \ox_A M$ and 
	hence on the cohomology $\Tor_A(M,N)$
	of the resulting single complex.
	It is the 
	total degree $n = p+q$
	that is the relevant grading 
	for \Cref{thm:Kapovitch} and \Cref{thm:Cartan}.
}
This can also be recovered through \eqref{fig:Kapovitch},
which is pullback square
because of the homeomorphism
$EG \ox_H G/K \lt EG/H \x_{EG/G} EG/K$
given by $e \ox gK \lmt (eH,egK)$.
If we take $\rho_H^*\: \H(BG) \to \H(BH)$
as a Sullivan model of $\rho_H$
and the left module structure map
$\H(BG) \inc R \ox_{\H(BG)} \H(BK)$
as a Sullivan model of $\rho_K$,
then the standard result on pullbacks 
again gives \Cref{thm:Kapovitch-detail}.

Both \eqref{fig:Cartan} and \eqref{fig:Kapovitch}
are pullbacks of fibrations \eqref{fig:genpull}.
The common features of our hypotheses
are that $\kk$ is a field of characteristic $0$
that $\H(\chi)$ itself is a Sullivan model for $\chi$,
and
that the Sullivan model $M_E$ for $E$ 
is an $\H(X)$-module resolution of $\H(E)$
in such a way that $\pi\: \H(X) \to M_E$ is the module structure map.
The common conclusion is a ring isomorphism $\Tor_{\H(B)}(\H X,\H E) \to \H(Y)$.
 
We would like something like this to hold more generally,
but in general, one cannot hope to find models with trivial differential,
or in other characteristics,
even spaces with polyonomial cohomology.
In general, one at least
has a map
$\Tor^0_{\H(B)}(\H X, \H E) =  \H(X) \ox_{\H(B)} \H(E) \lt \H(Y)$,
which is an isomorphism if $B$ is contractible and $\kk$ is a field,
by the K\"unneth theorem.
We still have cochain algebras
and a map $\C(X) \ox_{\C(B)} \C(E) \lt \C(Y)$
for all~$\kk$,
and this is again a quasi-isomorphism when $B$ is contractible,
by the Eilenberg--Zilber theorem,
but usually not otherwise.
To generalize the description of the cohomology of a pullback 
to cases with noncommutativity, 
nonzero differentials, and nonzero $B$,
Eilenberg and Moore consider a 
a refined notion of resolution.

\begin{definition}
Given a \DGA $A$, 
a \defm{\DG} \defd{$A$-module} $M$ 
is a \DG $\kk$-module $M$
which is simultaneously an $A$-module
in such a way that the 
action map $A \ox M \lt M$
is a cochain map (meaning $d(am) = da\.m + a\.(-1)^{|a|}dm$).
A map of \DG $A$-modules is an $A$-module map
which is also a degree-$0$ cochain map.
A sequence of 
$\DG$ $A$-module maps $P_p \to P_{p+1}$
is 
\defd{proper exact}
if for each fixed degree~$j$
the three sequences 
$P_\bul^{j}$,
$B^j(P_\bul)$,
and $H^j(P_\bul)$
of $\kk$-modules 
are exact.
A \DG $A$-module $P$ is \defd{proper projective}
if for each 
proper exact sequence $M \os {\smash f}\to N \to 0$ of $A$-modules,
each \DG $A$-module map $g\: P \lt  N$ lifts 
along~$f$.\footnote{\ 
	\emph{I.e.},
	there exists 
	a \DG $A$-module map $\wt g\: P \lt M$ with $f\wt g = g$.
	}
Given a \DG $A$-module $M$,
there always exists 
a \defd{proper projective resolution} of $M$,
a sequence of proper projective $\DG$ $A$-modules 
$P_\bul = (P_p)_{p \leq 0}$
and a \DG $A$-module map $P_0 \to M$
such that the extended sequence $P_\bul \to M \to 0$
is proper exact.
This guarantees 
$\H(P_\bul)$ be a projective 
$\H(A)$-module resolution of $\H(M)$.\footnote{\ 
	For $A = \kk$, this recovers
	the notion of a 
	\emph{Cartan--Eilenberg resolution} of a complex%
	~\cite[Ch.~XVII]{cartaneilenberg},
	used to define hypercohomology:
	given an additive functor $F$,
	one takes the single complex $P$ associated 
	to $P_\bul$
	and defines the hypercohomology of $M$ with respect to $F$
	to be $\H F(P)$.
	To relate this to later terminology,
	in one of the model structures
	on the category of half-plane bicomplexes $C^{\bul,\bul}$
	a cofibrant replacement of a single complex $C^{0,\bul}$
	is exactly a Cartan--Eilenberg 
	resolution~\cite[\SS4]{muroroitzheim2019},
	so that hypercohomology with respect to $F$
	is the derived functor of $F$.
}
\end{definition}

If $M$ and $N$ are \DG $A$-modules
and $P_\bul$ a proper projective resolution of $M$,
then as in the classical case,
$P_\bul \ox_A N$ inherits an internal grading $q$
defined on pure tensors by $q = |x \ox y| = |x| + |y|$,
the resolution grading $p$, 
and the total grading $n = p + q$.
The \defd{differential Tor}
$\defm{\Tor_A(M,N)}$,
defined as 
the cohomology of the associated single complex,
with the inherited bigrading, 
does not depend on the choice of $P_\bul$.

Returning to \Cref{fig:genpull},
a proper projective resolution $P_\bul$ 
of the $\C(B)$-module $\C(X)$
comes with a surjection $P_\bul \epi P_0 \epi \C(X)$
inducing a composite
\quation{\label{eq:EM-res}
	P_\bul \ox_{\C(B)} \C(E) \lt \C(X) \ox_{\C(B)} \C(E) \lt \C(Y)\mathrlap.
}

\bthm[Eilenberg--Moore]\label{thm:EM}
Suppose that 
$\kk$ is a principal ideal domain,
$\pi_1( B)$ acts trivially on the homotopy fiber $F$ of $E \lt B$,
and either 
(a) each of the $\kk$-modules $H^n(B)$ and $H^n(X)$ 
is finitely generated or 
(b) each $H^n(F)$ is.
Then \eqref{eq:EM-res} induces an isomorphism
\[
\Tor_{\C (B)}(\C X, \C E) \isoto \H(Y)\mathrlap,
\]
which is multiplicative under a certain natural ring structure on the domain.
\ethm
\bpf[Sketch]
Replace everything with a CW-complex
and
$X \lt B$ with a cellular map,
Filter $\C(Y)$ of \eqref{eq:EM-res} by 
the $X$-skeletal Serre filtration
$F_p \C(Y) = \ker\big(\C(Y) \to \C(\varpi^*X_{p-1})\big)$
and
$\C(E)$ by the $B$-skeletal Serre filtration, 
$P^\bul$ by resolution degree, and 
the domain of \eqref{eq:EM-res} by the tensor filtration.
Then \eqref{eq:EM-res}
is filtration-preserving
and hence induces a map of filtration spectral sequences,
whose codomain is the \SSS
converging to $\H(Y)$
and whose $E_2$ page map unwinds as
the identity map of $\H(X;\H F)$.
Thus
\eqref{eq:EM-res}
is a quasi-isomorphism.
\epf

If we instead filter
$P_\bul \ox_{\C(B)} \C(E)$ 
or more generally $P_\bul \ox_A N$
by the resolution degree $p$,
we get a left--half-plane spectral sequence of K\"unneth type,
the \defd{algebraic Eilenberg--Moore spectral sequence},
with $E_2 = \Tor_{\H (A)}(\H M, \H N)$, 
converging to $\Tor_{A}(M,N)$.

\bcor[\defd{Eilenberg--Moore spectral sequence} (\textcolor{RoyalBlue}{\EMSS})%
~{\cite{EMSS1965,smith1967emss}}]
Under the hypotheses of \Cref{thm:EM},
there exists a spectral sequence of $\H(B)$-algebras 
with 
$E_2 = \Tor_{\H (B)}(\H X, \H E)$,
converging to $\H(Y)$.%
\footnote{\
	The earliest version of differential homological algebra 
	and hence of the Eilenberg--Moore spectral sequence 
	was worked out by Eilenberg and Moore as early as 1957;
	per a 1959 lecture of Moore in the 
	\emph{S\'eminaire Henri Cartan}~\cite[fn.~1]{moore1959algebre},
	it was a topic of the Princeton topology seminar in 1957--8.
	The cohomological version 
	of differential homological algebra and the \EMSS
	appears in full in the 
	unpublished 1962 version 
	of Paul Baum's thesis~\cite{baumthesis}.
The first account of the cohomological \EMSS published in a journal
seems to have been Larry Smith's from 1967~\cite{smith1967emss}.
The homological version, with Cotor and coalgebras in place of Tor
and algebras, appears in the Eilenberg--Moore paper only in 1965~\cite{EMSS1965},
	the cohomological version being deferred to a Part~II yet to appear.
	Reference to the related bar spectral sequence
	appears in work of Clark~\cite{clark1965homotopy} to be discussed later.
}
\ecor

There is typically no easy way to compute the differentials 
of the \EMSS unless it is known to collapse,
so many authors have set themselves the task of proving \EMSS collapse results.
A recurrent strategy runs through the algebraic \EMSS:

\bprop[{%
\cite[XI.3.2]{maclane}%
\cite[Cor.~1.8]{gugenheimmay}%
\cite[Theorem~5.4]{munkholm1974emss}}%
]%
\label{thm:EMSS-iso}
The algebraic \EMSS associated to a diagram
$M \longfrom A \lt N$ of maps of nonnegatively-graded \DGAs 
is convergent and functorial 
in the sense that
a commutative diagram

\vspace{-2ex}

\quation{\label{eq:Tor-DGA-functoriality-squares}
	\begin{aligned}
		\xymatrix@C=2.25em{
			M' \ar[d]_(.45)u
			&	\ar[l]_{\phi_{M'}} \ar[r]^{\phi_{N'}}
			A'
			\ar[d]|(.375)\hole|(.45)f|(.525)\hole
			& 	N\mathrlap'
			\ar[d]^(.45)v\\
			M
			& A
			\ar[l]^{\phi_M} \ar[r]_{\phi_N}	
			& N
		}
	\end{aligned}
}
of \DGA maps
induces a map of spectral sequences.
In particular, 
	if $f$, $u$, $v$ are quasi-isomorphisms,
the induced map	
$	\Tor_f(u,v)$
is a graded-linear isomorphism.
\eprop
\nd 
Thus if we can find vertical maps 
making the diagram 
	\begin{equation}\label{eq:Munkholm-squares-sketch}
	\begin{aligned}
		\xymatrix@C=1em@R=2.5em{
			\H (X)
				\ar[d]
			&	
			\H (B)
				\ar[d]
				\ar[l]
				\ar[r]
			&
			\H (E) 
				\ar[d]
			\\
			\C(X)
			&
			\C(B)
				\ar[l]
				\ar[r]
			&
			\C(E)
			\mathrlap,
				}	
	\end{aligned}
	\end{equation}
commute,
we can conclude $\H(Y) \iso \Tor_{\C(B)}(\C X,\C E) \iso 
\Tor_{\H(B)}(\H X,\H E) $,
a strong \EMSS collapse result,
and in the case that \Cref{fig:genpull}
is \Cref{fig:Cartan} or \Cref{fig:Kapovitch},
we obtain \Cref{thm:Cartan} or \Cref{thm:Kapovitch}, respectively.
We usually 
cannot show \eqref{eq:Munkholm-squares-sketch}
is commutative on the nose, but
we will encounter more general versions of Tor,
each admitting a functorial algebraic \EMSS
and an analogue of \Cref{thm:EMSS-iso},
and for these generalizations,
we will be able to follow this strategy.

\section{Baum}\label{sec:Baum}
Paul Baum's 1962 thesis~\cite{baumthesis} 
aimed to establish \Cref{thm:Cartan} for $\kk$ a field
by proving the collapse
of the \EMSS of \Cref{fig:Cartan}.
For $T$ a maximal torus of $K$,
he noted
there is a map of fibrations from $G/T \to BT \to BG$
to $G/K \to BK \to BG$,
inducing a map of {\EMSS}s
which he showed was of the form $E_r \inc E_r \ox \H(K/T)$.
As a consequence,
one has the following.

\bthm[Baum {\cite[3.3.2]{baumthesis}}]\label{thm:Baum-lemma}
Let $\kk$ be a field. 
If $T$ is a maximal torus of $K$,
the \EMSS converging to $\H(G/K)$ collapses if and only if
that converging to $\H(G/T)$ does.
\ethm
\nd Later proofs often call on mild variants of \Cref{thm:Baum-lemma};
for the purposes of this survey, we will gloss all of them 
as \defd{Baum's reduction}.

Unfortunately, 
there is an error in the proof of the 
main collapse result.\footnote{\ 
	Paul will readily tell you this if 
	you happen to coincidentally sit next to him at a conference 
	and tell him you've just read his thesis.
}
In the 1968 published version~\cite{baum1968homogeneous},
the following is salvaged.

\bthm[Baum]
For any field $\kk$,
the \EMSS associated to \Cref{fig:Cartan} collapses
when the kernel of the map 
$Q\H(BG) \to Q\H(BT)$ of indecomposables
has dimension $\leq 2$.
In particular, \Cref{thm:Cartan} holds additively.
\ethm
\bpf
The \EMSS 
is concentrated in even rows $q$,
forcing $d_2 = 0$.
By the assumption on indecomposables,
it is generated in columns $-2 \leq p \leq 0$,
forcing $d_{\geq 3} = 0$.
\epf
In particular, when $G$ and $K$ are of equal rank,
	$\H(BG) \lt \H(BK)$ is surjective 
	so the \EMSS is concentrated in 
	the $0\th$ column,
	forcing a \CGA isomorphism $\H(G/K) \isoto
	\kk \ox_{\H(BG)} \H(BK)$
recovering Borel's result~\cite[Cor.~7.5]{baum1968homogeneous}.


\bex[Borel]
Consider $K = \SU(5) < \U(5) < \Sp(5) = G$.
Computing the map $\H B\Sp(5) \lt \H B\SU(5)$
and computing Tor using a Koszul complex,
one has $\H(G/K) = \F_p\big\{1,c_3, c_5,w_{21},w_{26},[G/K]\big\}$,
where $c_j \in H^{2j}(G/K)$ are the images of the universal 
Chern classes, $\dim G/K = 31$,
and the only nonzero products are those implied by Poincar\'e duality,
for $\kk = \F_p$ ($p \neq 2$) or $\kk = \Q$.
This looks different than the previous results because
the ideal of $\H B\SU(5)$ generated by $H^{\geq 1}B\Sp(5)$
is not generated by a regular sequence,
and equivalently (for $\kk = \Q$),
$G/K$ is not formal in the sense of rational homotopy theory.
For $\kk = \F_2$,
one instead finds $\H(G/K) \iso \Lambda[z_3] 
\ox \F_2[c_2,c_3,c_4,c_5]/(c_2^2,c_3^2,c_4^2,c_5^2)$.
\eex

\begin{counterexample}\label{thm:Baum-cex}
The center of the unitary group $\U(2)$ 
is the diagonal copy $\D\U(1)$ of $\U(1)$,
which meets $\SU(2)$ in $\pm I$.
We thus have diffeomorphisms 
$\SO(3) \iso \SU(2)/\{\pm I\} \iso \U(2)/\D\U(1)$,
and we consider the \EMSS of 
$\SO(3) \to B\D\U(1) \to B\U(2)$ over $\F_2$
beginning with \[E_2 = \Tor_{\H B\U(2)}\big(\F_2,\H\D\U(1)\big)\mathrlap.\]
To compute this, 
one can check 
the differential on $\H B\D\U(1) \ox \U(2) = 
\F_2[y_2] \ox \Lambda[z_1,w_3]$
takes~$z_1$ to~$0$ and $w_3$ to $y_2^2$,
so the $\Tor$
is isomorphic to $\Lambda[z_1] \ox \F_2[y_2]/(y_2^2)$.
This is isomorphic to
$\H \SO(3) = \F_2[x_1]/(x_1^4)$
as a graded vector space but not as a ring,
so the full multiplicative version of \Cref{thm:Cartan}
is not true for $\kk = \F_2$.
\end{counterexample}

\section{Cup-\texorpdfstring{$i$}{i} products}\label{sec:cup-i}
All later results showing \EMSS collapse
require notions of comparisons of $\H(X)$ and $\C(X)$
by linear maps that are multiplicative only up to homotopy.
If $\kk[x,y]$ is the cohomology ring of some space $X$,
finding an isomorphic copy of $\kk[x,y]$ 
in $\C(X)$ itself is generally impossible,
because
lifting to representatives $x,y \in \C(X)$,
one only has $x \cup y \equiv (-1)^{|x||y|} y \cup x$ 
modulo a coboundary.
This coboundary can be chosen 
in such a way
as to yield a cochain homotopy from the cup product to  
$(x,y) \lmt (-1)^{|x||y|} y \cup x$,
 called the Steenrod \defd{cup-$1$ product} 
and denoted~$\defm{\cupone}$%
~\cite[\SS\SS2,\,5]{steenrod1947products}.
The cup-$1$ product 
is itself commutative up to a cochain homotopy
witnessed by an operation $\cuptwo$, 
and inductively Steenrod found 
a sequence of binary operations ${\cup}_{ i}$ of degree $-i$, 
each commutative up to a homotopy witnessed by ${\cup}_{ i+1}$.

\begin{definition}\label{def:cup-i}
	Given \DGAs $A$ and $B$,
	there is a natural \DGA isomorphism
	$\chi_{A,B}\:$ $A \ox B \lt B \ox A$ 
	given on pure homogeneous tensors
	by $a \ox b \lmt (-1)^{|a||b|} b \ox a$.
	A \DGA $A$ is said to \defd{admit cup-}$\defm i$ \defd{products}
	if for $0 \leq j \leq i$
	there exist degree-($-j$) operations $\mu_j\: A \ox A \lt A$,
	starting with $\mu_0 = \mu_A$ the ring multiplication of $A$,
	such that for $0 \leq j < i$ 
	one has
	$
	D\mu_{j+1} = \mu_j - \mu_j\chi_{A,A}
	\mathrlap.
	$
\end{definition}

\section{May}\label{sec:May}

Peter May characterized 
the differentials of the algebraic \EMSS in terms of generalized 
Massey products defined using matrices of cochains, 
called \emph{matric Massey products},
and showed that in the \EMSS associated to a bundle $F \to E \to B$,
the elements of $\C(B)$ figuring in these matric Massey products
were iterated $\cupone$-products of cocycles.
He announced~\cite{may1968principal} 
that he had found a \DGA quasi-isomorphism
$f\: \C(BT) \lt \H(BT)$ for any $\kk$,
which by the commutativity of $\H(BT)$ annihilates $\cupone$-products.%
	\footnote{\ 
	He also found
	\DGA quasi-isomorphisms $\C\big(\mn\prod K(\pi_j,n_j);\F_2\big) 
	\lt \H\big(\mn\prod K(\pi_j,n_j);\F_2\big)$
	for positive integers $n_j$
	and finitely generated abelian groups $\pi_j$
	(with no $4$-torsion for $n_j = 1$).
}
This induces a map $\Tor_{\id}(\id,f)$
from $\Tor_{\C(BG)}(\kk,\C BT)$ 
to $\Tor_{\C(BG)}(\kk, \H BT)$,
which is an isomorphism by \Cref{thm:EMSS-iso}.
The differentials in the \EMSS converging to $\Tor_{\C(BG)}(\kk, 
\textcolor{red}{\C} BT)$
all involve $\cupone$-products on $\C(BT)$,
which $f$ annihilates, 
showing the algebraic \EMSS converging to $\Tor_{\C(BG)}(\kk, \textcolor{red}{\H}BT)$
collapses at $E_2 = \Tor_{\H(BG)}(\kk, \H BT)$.
With Baum's reduction,
this establishes 
the following:
\bthm
Let $\kk$ be Noetherian.
Then \Cref{thm:Cartan} holds
additively up to an extension problem;
\emph{i.e.}, with the filtration on $\H(Y) \iso 
\H\big(P_\bul \ox_{\C(BG)} \C(BK)\big)$ induced by the resolution degree~$p$,
the associated graded module $\gr \H(Y)$ is isomorphic
to $\Tor_{\H(BG)}(\kk,\H BK)$.
\ethm
A preprint was circulated,
but the proofs were involved 
and did not see print,
suppressed in favor of later proofs
joint with V.K.A.M Gugenheim.

\section{Gugenheim--May}\label{sec:Gugenheim-May}

Gugenheim and May~\cite{gugenheimmay} construct a \DGA formality map 
$f\: \C(BT) \lt \H(BT)$ annihilating $\cupone$-products
by dualizing \DG Hopf algebra quasi-isomorphisms
from homology $H_*(BT)$ to 
$C_*(BT)$, 
where the particular model of 
$BT$ is a direct power of a simplicial model of $BS^1$.
They assume only that $\kk$ is a Noetherian ring.
Then,
using a variant of Baum's reduction \ref{thm:Baum-lemma}
requiring a torsion hypothesis on~$BK$,
they reduce the \EMSS collapse 
for \Cref{fig:Cartan}
to that 
for $G_T \to BT \os\rho\to BG$.

As May had earlier observed,
the map $\Tor_{\id}(\id,f)\:
\Tor_{\C(BG)}(\kk,\C BT) \to
	\Tor_{\C(BG)}(\kk, \H BT)$ 
is a linear isomorphism.
To compute the codomain,
Gugenheim--May resolve $\kk$ as a $\C(BG)$-module in the following way.
Fix exterior generators $z_j \in \H(G)$,
corresponding polynomial generators 
	$x_j = \tau z_j\in \H(BG)$,
	and representatives
	$c_j \in \C(BG)$.
Equip the bigraded module
	$\defm M = \C(BG) \ox \H(G)$
	with a differential $\defm{d_M}$
	whose value on each pure tensor
	$1 \ox z_{j_1} \cdots z_{j_n}$
 	differs from the naive choice
	$\sum_i \pm c_{j_i} \otimes z_{j_1} \cdots \wh{ z_{j_i}} \cdots z_{j_n}$ 
	(the caret~$\wh{\ }$ denoting omission)
	by an element of the ideal 
	generated by $\cupone$-products on $\C(BG)$.\footnote{\ 
This is the best one can do, since $\C(BG)$ is not actually commutative.
	}
	Resolving $\kk$ by $M$,
	one can compute
	$\Tor_{\C(BG)}(\kk,\H BT)$
	from 
	$
	\H(BT) \ox_{\C(BG)} M$
	$
	\iso
	\H(BT) \ox \H(G)$,
	where the
	differential annihilates $\H(BT)$
	and takes $1 \otimes \prod z_{j_i}$ to
	$(f\mn \rho^* \otimes \id) d_M(1 \otimes \prod z_{j_i})$.
	But $f\mn \rho^*$ annihilates the $\cupone$-products
	distinguishing these values from those
	of the Cartan/Borel model,
	defined by $1 \ox z_j \lt \rho^*x_j \ox 1$,
	so we have recovered the classical $\Tor_{\H(BG)}(\kk,\H BT)$.
This direct computation does not pass through an \EMSS collapse
result and thus also resolves the additive extension problem.

\bthm[Gugenheim--May]\label{thm:GuMa}
Let $\kk$ be a Noetherian ring 
such that $\H(BK;\Z)$ has no $p$-torsion 
for any factor $p$ of $\chara\kk$ (but is not necessarily polynomial).
Then \Cref{thm:Cartan} holds additively.
\ethm

\section{A-infinity notions historically}\label{sec:A-infty}

We have seen the assumption of a cup-one product on a \DGA $A$ 
limits how badly a map $\kk[x,y] \lt A$ lifting generators
can fail to be a ring map.
It is not hard to check a \DGA
$A$ is commutative
if and only if the multiplication
$\mu\: A \ox A \lt A$, a cochain map,
is in fact a \DGA map,
so a more systematic way of limiting noncommutativity 
is a system of homotopies
moderating $\mu$'s failure to be multiplicative.

Notions of \emph{strong homotopy multiplicativity} (\textcolor{RoyalBlue}{\SHM})
and \emph{strong homotopy commutativity} 
(\textcolor{RoyalBlue}{\SHC})
originally apply to topological monoids
and are due to Sugawara~\cite{sugawara1960homotopy}.
We cannot discuss them in detail here,
but
%
Sugawara was able to show that an 
\SHM map between two monoids $G$ and $H$ 
induces a map $BG \lt BH$ of classifying spaces,
that a countable CW-complex $B$
has loop space~$\Omega B$ \SHC if and only if $B$ is an H-space,
and that for $G$ a topological group,
$BG$ is an H-space if and only if $G$ is \SHC.
Allan Clark algebraized these notions~\cite{clark1965homotopy}\footnote{\,$\phantom{j}\!\!\! $
He is interested 
in the monoid $\Omega X$
of variable-length (\emph{Moore}) loops 
$\coprod_{r \geq 0} \Map\Big(\big([0,r],\{0,r\}\big),(X,*)\Big))$
on a pointed space $X$
and particularly in comparing the monoids 
$\Omega (X \x X)$ and $\Omega X \x \Omega X$.
These are homotopy equivalent but not homeomorphic,
which had lead to an error 
in a lecture of Moore~\cite[Thm.~7.II, pf.]{moore1959algebre}.
}
weakening Sugawara's hypotheses but obtaining
similar consequences.
He defines an \emph{\SHM map} of \DGAs $A \lt A'$
as what is now called an \defd{\Ai-algebra map},
a sequence of maps $A^{\otimes n} \lt A'$
satisfying certain homotopy coherence conditions
approximating multiplicativity.
He shows these amount to a map of differential
graded coalgebras  $\B A \lt \B A'$,
where $\B A$ is the \emph{bar construction}
defined in the next section.

Stasheff introduced associahedra $K_n$
and $A_n$-spaces,
and showed an $A_n$-space is equivalently a space $X$ 
endowed with a sequence of 
maps $K_j \x X^{j} \lt X$ for $j \in [2,n]$ 
satisfying now--natural-seeming conditions~\cite{stasheff1963I,stasheff1963II}.
He defined an $A_n$-map in such a way that when $n = \infty$,
it is an \SHM map in the sense of Sugawara.
An \defd{$A_n$-algebra} is an augmented graded module equipped
with a sequence of linear maps $m_j\: A^{\ox j} \lt A$ for $j \in [1,n]$
satisfying formally similar conditions.
	Stasheff showed that an $A_n$-algebra is, 
	in later language,
	a module for the $n\th$ filtrand 
	of the operad\footnote{\ 
		We do not assume or use any specific results
		about operads in this survey,
		so we state this merely for historical
		context and those who already know about operads.
	}
	of cellular chains on $K_\bul$,
	so that if $X$ is an $A_n$-space, 
	$C_*(X)$ is an $A_n$-algebra.
An \defd{\Ai-algebra} is an augmented chain complex $A$ 
which is an $A_n$-algebra for all $n$,
amounting to a differential making the tensor coalgebra 
$\Direct_{p=0}^\infty {\ol A}{}^{\otimes p}$
a differential graded coalgebra.
When $A$ is a \DGA,
this prescription
gives the differential on the bar construction 
of the next section,
up to sign.
Stasheff also connects the algebraic and topological 
bar constructions via
a chain equivalence $\B C_*(G) \lt C_*(BG)$,
for $G$ a topological monoid,
which will be used in \Cref{sec:HSM}.

\section{The bar construction}\label{sec:bar}

The bar construction is a functor from $\Algs$
to a category $\Coalgs$
of cochain complexes satisfying axioms dual to the axioms for a ring 
with a derivation.
%
\begin{definition}
	A \defd{coaugmented differential graded $\kk$-coalgebra}
(\textcolor{RoyalBlue}{\DGC})
is a cochain complex $(C,d_C)$
equipped with a \defd{comultiplication} $\defm{\D_C}\: C \lt C \ox C$
and maps $\kk \os{\h_C}\to C \os{\e_C}\to \kk$ composing to $\id_\kk$,
respectively the \defd{coaugmentation} and \defd{counit},
satisfying the identities
\[
(\id \ox \D)\D = (\D \ox \id)\D,\quad\ 
\D d = (d \ox \id + \id \ox d)\D,\quad\ 
(\e \ox \id)\D = \id = (\id \ox \e)\D\mathrlap.
\]
A (coaugmentation-preserving) \defm{\DGC} \defd{map} is a degree-$0$
cochain map $g\: C \lt C'$ 
satisfying $\D_{C'} g = (g \otimes g)\D_C$
and $\e_C' g = \e_C$ (and $g\h_C = \h_C'$).
One writes $\D c = \sum \defm{c_{(1)}} \ox \defm{c_{(2)}}$.
A tensor product $C \ox C'$ of \DGCs 
becomes a \DGC under the tensor differential
and  comultiplication
taking $c \ox c'$ to $\sum (-1)^{\smash{|c_{(2)}||c'_{(1)}|}}
c_{(1)} \ox c'_{(1)} \ox c_{(2)} \ox c'_{(2)}$.
When $\D$ is itself a \DGC homomorphism $C \lt C \ox C$,
then $C$ is called \defd{cocommutative}.
Write 
$\D^{[n+1]} =
(\D \otimes \id^{\otimes n-1})\D^{[n]}\: C \lt C^{\otimes n}$
for the iterates of 
$\defm{\D^{[2]}} = \D_C$
and
$\smash{\defm{\ol C}} = \coker \h_C$
for the coaugmentation coideal;
a \DGC $C$ is \defd{cocomplete}
if each homogeneous element 
is annihilated by one of the composites
$C \to C^{\otimes n} \to \ol C{ }^{\otimes n}$.
Write
$\defm{\Coalgs}$ for the category of cocomplete \DGCs, 
and coaugmentation-preserving \DGC maps.
%
A \defm{\DG} \defd{Hopf algebra}
is a \DGA $A$ equipped with a \DGA homomorphism $\D\: A \lt A \ox A$
(equivalently, a \DGC $A$ with a \DGC map $\mu\: A \ox A \lt A$).
\end{definition}

\bex
The singular chain complex $C_*(X)$ associated to a topological space $X$
becomes a \DGC under the map 
taking a singular simplex $\s\: \D^n \to X$
to the sum of $\s|_{\D^{[0,\ldots,p]}} \ox \s|_{\D^{[p,\ldots,n]}} 
\in C_p(X) \ox C_{n-p}(X)$ for $0 \leq p \leq n$.\footnote{\ 
	This is the case $Y_\bul = \Sing\, X$ 
	of a more general 
	\DGC structure
	on the chain complex $C_*(Y_\bul)$ associated to a simplicial
	set $Y_\bul$.}

Using the Eilenberg--Zilber \DGC quasi-isomorphism
if $X$ is an H-space, the 
composite \DGC map 
$C_*(\mu_X) \o \nabla$
makes $C_*(X)$ a \DG Hopf algebra, 
and under sufficient flatness hypotheses,
$H_*(X)$ becomes a cocommutative Hopf algebra.
Dually, $\C(X)$ becomes a \DG Hopf algebra 
and $\H(X)$ a commutative Hopf algebra.
\eex

%

\begin{definition}\label{def:BA-structure}
The \defd{desuspension} $\defm{\desusp \ol A}$ of 
the submodule $\defm {\ol A} = \ker \e_A$
of an augmented cochain complex $A$
is $\ol A$ regraded via $(\desusp \ol A)_n \ceq \ol A_{n+1}$,
given differential $d_{\desusp\ol A} =  -\desusp d_A \susp$.
On the direct sum $\defm{\B A}$ of $\defm{\B_p A} \ceq
(\desusp \ol A)^{\otimes p}$ ($p \geq 0$),
the \defd{tensor coalgebra} structure
takes 
$\desusp a_1 \otimes \cdots \otimes \desusp a_p \eqc
\defm{[a_1|\cdots|a_p]} \in \B_p A$
to 
$\sum_{0 \leq \ell \leq p} 
[a_1|\cdots|a_\ell] \otimes [a_{\ell+1}|\cdots|a_p]
\in \B A \ox \B A$,
where $[] = 1 \in \kk = \B_0 A$.
When~$A$ is a \DGA,
the \defd{bar construction} 
is the \DGC structure on $\defm{\B A}$ 
whose differential~$d_{\B A}$
is the sum of the tensor differentials 
on the $\B_p A = (\desusp \ol A)^{\otimes p}$
and the ``bar-deletion'' maps 
$\id^{\otimes i} \otimes \desusp \mu (\susp \otimes \susp) \otimes \id^{\otimes j}$
on $\B_{i+j+2}$ ($i,j\geq 0$)
taking $[a|b|c] \lmt \pm [a|bc] \pm [ab|c]$ and so on.%
\footnote{\ 
	Signs are all determined by the Koszul convention
	$(f \ox g)(a \ox b) = (-1)^{|g||a|} f(a) \ox g(b)$.
}
\end{definition}
	
%

\nd We can use the bar construction to generalize the notion 
of commutativity.
As we have noted,
a \DGA
$A$ is commutative
if and only if $\mu\: A \ox A \lt A$ is itself a \DGA map.
Weakening this requirement,
Stasheff--Halperin~\cite[Def.~8]{halperinstasheff1970}
call a \DGA $A$
an \textcolor{RoyalBlue}{\SHC}\defd{-algebra}
when it admits
a \DGC map $\Phi\:\B(A \ox A) \lt \B A$
such that $\Phi[a \ox b] = [ab]$ for $a \ox b \in \ol{A \ox A}$---%
asking, in other words, 
only that~$\mu$
be the unary component $(\ol{A\otimes A})^{\otimes 1} \lt A$ 
of an \Ai-algebra map.\footnote{\ 
    Clark's definition~\cite{clark1965homotopy} 
    had asked even less,
    merely that  
    $\Phi[a\ox 1] = [a] = \Phi[1 \ox a]$.
}


\begin{definition}[See Husemoller \emph{et al.}~{\cite[Def.~IV.5.3]{husemollermoorestasheff1974}}]%
\label{def:shuffle}
There exists a natural transformation
\[
		\defm\nabla \: \B A_1 \ox \B A_2 \lt \B(A_1 \ox A_2)
\]
\nd of functors $\Algs \x \Algs \lt \Coalgs$,
the 
\defd{shuffle map}, 
which is a homotopy equivalence of cochain complexes.
	It is the direct sum of the maps $\BBB_p A \ox \BBB_q B \lt \BBB_{p+q}(A \ox B)$
	sending $[a_1|\cdots|a_p] \otimes [b_1|\cdots|b_q]$ 
	to the sum of all tensor $(p,q)$-shuffles (with Koszul sign) of
	$[a_1 \otimes 1|\cdots|a_p \otimes 1|1 \otimes b_1|\cdots|1 \otimes b_q]$.
\end{definition}

\nd Clark notes that the composite $\Phi \o \nabla$
makes $\B A$ a (possibly nonassociative) \DG Hopf algebra.
The earlier Eilenberg--Mac Lane paper~\cite[(7.7)]{eilenbergmaclane1953}
already showed this for $\Phi = \B\mu$
when $A$ is a  \CDGA. 
We will strengthen this observation in \Cref{thm:Franz-SHC}.

The bar construction $\B A$ for a \CDGA $A$ was introduced in print
by Eilenberg--Mac Lane~\cite[\SS11]{eilenbergmaclane1953}%
\footnote{\ 
    This is the same work where
    the Eilenberg--Zilber map $\nabla$ is introduced;
    the Eilenberg--Zilber theorem was originally proved
    using acyclic models.
    The bar shuffle $\nabla\:\B A_1 \ox \B A_2 \lt \B(A_1 \ox A_2)$
    and the observation it is a homotopy equivalence
    first appeared in Eilenberg--Mac Lane's sequel%
    ~\cite[Thm.~4.1a]{eilenbergmaclane1954}.
}
not to parameterize homotopy-associative operations,
but to provide a functorial resolution of a \DG $A$-module $M$,
proper projective  
when $A$, $M$, and $\H(A)$
are all flat over $\kk$.

\bdefn\label{def:one-sided}
The \defd{one-sided bar construction} of
a \DGA $A$ and a right \DG $A$-module~$M$
is the graded $\kk$-module $M \otimes \BA$
(note $M \otimes \B_0 A = M \otimes\,  \kk = M$)
equipped with  
the sum of the tensor differential
and the maps 
$\defm{\ell_p}\: \mu_M(\id_M \otimes s) \otimes
\id^{\otimes p -1} \: M \otimes \B_p A  \lt M \otimes \B_{p-1} A$
taking $m[a_1|\ldots|a_{p-1}|a_p] \lmt 
\pm m a_1[a_2|\ldots|a_{p}]$.
Given a right $A$-module $N$,
the \defd{two-sided bar construction}
is $(M \otimes \B A) \otimes_A N$,
which can be identified with $M \otimes \B A \otimes N$
with the sum of the tensor differentials,
the $\ell_p \ox \id_N$,
and the maps
$\defm{r_p} \: 
-\id_M \otimes \id^{\otimes p -1} \otimes \mu_N(s \ox \id_N)
\: M \otimes \B_p A \otimes N \lt M \otimes \B_{p-1} A \otimes N$
taking $m[a_1|\cdots|a_p]n \lmt \pm m[a_1|\cdots|a_{p-1}]a_p n$.
\edefn

Mild flatness hypotheses imply 
$\Tor_A(M,N) = \H \B(M,A,N)$.

\section{Stasheff--Halperin}\label{sec:SH}
%
%
Stasheff--Halperin\footnote{\ 
	written in this order; 
	the pun in the title is also deliberate
	and similarly unexplained}~\cite{halperinstasheff1970},
in an offering to a workshop proceedings volume,
suggested a program
to generate collapse results for the \EMSS 
of \Cref{fig:Cartan}
using an \Ai-map 
they construct.
Given an \SHCA $(A,\Phi)$,
define the \defd{iterates} 
$\defm{\iter\Phi{n}}\: \B(\mnn A^{\otimes n}) \lt \B A$
of its structure map
by $\iter\Phi{2} \ceq \Phi$ and $\iter\Phi{n+1} \ceq \Phi(\iter\Phi{n} \T \id\mn_A)$.\footnote{\ 
    The iterates, 
    including $\Phi$, 
    are themselves \SHCA maps~\cite[Prop.~4.5]{munkholm1974emss}.
    }
From a list 
$(f_j\: A_j \to A)_{0 \leq j < \kappa}$ 
of \DGA maps,
one can define a composite
\DGC map 
\quation{\label{eq:compilation}
	\B\big(\mn\mnn\Tensor A_j\big) 
		\xtoo[\ \B(\Tensor  f_j)\!]{} 
	\B(\mnn A^{{\otimes} \kappa}) 
		\xtoo[\iter \Phi \kappa]{}
	\B A
}
combining the~$f_j$
when $\kappa$ is finite~\cite[Thm.~9]{halperinstasheff1970} 
(or,
via a careful colimiting 
argument,
countably infinite~{\cite[Prop.~3.9(i); \SS4.3]{munkholm1974emss}}). 
We will call this map the \defd{compilation} or \emph{assembly}
of the $\B f_j$.

When $A$ is an \SHCA with countably generated polynomial cohomology 
\[\H\mn A = \kk[x_1,x_2,\ldots] \iso \Tensor_{\smash j} \kk[x_j]\mathrlap,\]
a choice of representative $a_j \in A$ 
for each~$x_j$
induces a \DGA map $f_j\: \kk[x_j] \lt A$ sending $x_j$ to~$a_j$,
so \eqref{eq:compilation}
gives a \DGC map 
\quation{\label{eq:lambda}
	\defm{\l_A}\: \B\H\mn A 	\lt 	\B A\mathrlap.
}
The map $\l_A$ induces
the identity map $\H\mn A = \H(\H\mn A) \lt \H\mn A$ as follows.

\begin{proposition}[{\cite[Prop.~3.7]{munkholm1974emss}}]\label{def:H-flat}
	If $A$ and $B$ are \DGAs and $g\: \B A \to \B B$ a \DGC map,
	 the composite
     ${\ol A \us{\desusp}\lt \B A \us g\to \B B \epi \ol B}$
	induces a graded algebra map 
	\[
		\defm{\Hflat}g\: \H\mn A \lt \H B
	\]
such that $\Hflat \l_A = \id_{\H\mn A}$.
\end{proposition}

Stasheff--Halperin  describe an \SHM \emph{module}
over a \DGA $A$ as a \DG module $M$ equipped with 
linear maps $A^{\otimes j} \ox M \lt M$
whose adjoints $A^{\otimes j} \lt \End M$ 
yield an \Ai-map $\B A$ from $\B(\End M)$.
Using this data they define a certain differential on $\B A \ox M$
generalizing the one-sided bar construction of \Cref{def:one-sided}
and reuse the notation $\defm{\Tor_{A}(\kk,M)}$ 
for its cohomology.
Then for $\l_K$ and $\l_G$ the
$\Hflat$-quasi-isomorphisms of the previous paragraph
and $\defm\rho = B(K \inc G)$,
the composites
\begin{equation}\label{eq:two-composites}
	\begin{gathered}
	\xymatrix@R=-0.5em@C=5em{
&\B \H(BK)\ar[rd]^{\l_K}&\\
\B \H(BG) 	\ar[rd]_{\l_G}\ar[ru]^{\B\H(\rho)}&&\B \C(BK)\\
&\B \C(BG) \ar[ru]_{\B\C(\rho)}&
	}
	\end{gathered}
\end{equation}
induce two \SHM $\H(BG)$-structures on $\C(BK)$.
The maps $\l_K$ and $\l_G$ respectively
induce graded isomorphisms
\eqn{
\Tor_{\H BG}(\kk,\H BK) \isoto\ & \Tor^{\mr{top}}_{\H BG}(\kk,\C BK)\mathrlap,
\\
&\Tor^{\mr{bottom}}_{\H BG}(\kk,\C BK) \isoto \Tor_{\C BG}(\kk,\C BK)\mathrlap.
}
If the top and bottom \SHM $\H(BG)$-module
structures on $\C(BK)$ agreed, 
the composite
isomorphism
would imply \Cref{thm:Cartan} additively,
but they do not,
and later work (1) finds a homotopy between them
or (2) postcomposes a map annihilating the difference.
Both approaches 
involve the notion of a twisting cochain.

\section{Twisting cochains}\label{sec:history-cochains}

Twisting cochains 
were originally defined
by Edgar Brown~\cite{brown1959twisted},
who asked what algebraic information was necessary to generalize
the Eilenberg--Zilber quasi-isomorphism
\[\defm\nabla\: C_*(B) \ox C_*(F) \lt C_*(B \x F)\]
from a product 
to a fiber bundle
$F \to E \to B$.
He found that $C_*(E)$ 
was quasi-isomorphic to 
the chain complex given by a new 
``twisted'' differential on $C_*(B) \ox C_*(F)$
defined using the information
encoded in an inhomogeneous cochain
in $C^*(B;C_* \Omega B)$,
viewed as a map $t_B\: C_*(B) \lt C_{*-1}(\Omega B)$,
and the continuous action $\w\:\Omega B \x F \lt F$
via path-lifting.
This new differential and the differential of the one-sided
bar construction $\B A \ox M$ are both instances of a single construction
we will describe in \Cref{def:twisted-tensor-product}.


%
\begin{definition}[Brown~{\cite[\SS3]{brown1959twisted}}]%
For $C$ a \DGC and $A$ a \DGA,
the \defd{cup product}
$
\defm{f \cup g} \ceq \muA(f \tensor g)\D_C$
renders $\GM(C,A)$ a \DGA~{\cite[\SS1.8]{munkholm1974emss}}.
\end{definition}


%

\begin{proposition}\label{thm:bar-cofree}
Let $A$ be a \DGA.
The composite $\defm{t^A}\:
\B A \epi \B_1 A = \desusp \ol A \simto \ol A \inc A$
is a cochain map.
For any \DGC~$C$,
	\DGC maps $F\:C \lt \BBB A$ 
	correspond bijectively
	via $F \lmt t^A \o F$
to maps
$t \in \GM_1(C,A)$
satisfying the three conditions
\quation{\label{eq:twisting-cochain}
\eA t = 0 = t \h_C,\qquad \Homd t = t \cup t\mathrlap.
}
\end{proposition}

\begin{definition}[Brown~{\cite[\SS3]{brown1959twisted}}]%
	\label{def:twisting-cochain}
Let $C$ be a \DGC and $A$ a \DGA.
An element of $\GM_1(C,A)$ satisfying the conditions 
\eqref{eq:twisting-cochain} is a \defd{twisting cochain}.
We write $\defm{\Tw(C,A)}$ for the set of these.
The twisting cochain $t^A\: \B A \lt A$,
called the \defd{tautological twisting cochain},
is natural in the \DGA $A$.
\end{definition}

\begin{definition}[Brown~{\cite[\SS3]{brown1959twisted}}]%
	\label{def:twisted-tensor-product}
Let $C$ be a \DGC, 
$A$ a \DGA, 
$Q$ a differential right $C$-comodule,
and $N$ a differential left $A$-module.
We define the \defd{cap product} 
with an element~$\phi \in \GM(C,A)$ by
\eqn{
	\defm{\delta^{\mr R}_{\phi}} 
		\ceq 
	(\id_{Q}\otimes\mu_N)\,
	(\id_{Q}\otimes\, \phi \otimes \id_{N})\,
	(\D_Q \otimes \id_{N})\:
	Q \ox N &\lt Q \ox N,\\
	x \otimes y 
		&\lmt 
	\sum \pm x_{(1)} \otimes \phi(x_{(2)}\mnn)y.%
}
When $\phi = t$ is a twisting cochain,
$d_{Q} \ox \id_N + \id_Q \ox d_N  - \d^{\mr R}_{t}$
is a differential on 
$Q \ox N$. 
The resulting cochain complex
is 
the \defd{twisted tensor product} $\defm{Q \ox\twist{t} N}$.
Given another twisting cochain $t'\: C \lt A'$,
there is a similar ``left'' cap product $\defm{\d_{t'}^{\mr L}}$, 
for $M$ a right \DG module over a \DGA $A'$,
and one can define a twisted tensor product 
$M \ox_{t'} Q$ with differential differing from the tensor
differential by $\d_{t'}^{\mr L}$.
Taking $Q = C$, one can also form a two-sided twisted tensor product
$M \ox_{t'} C \ox_t N$.
\end{definition}

\bex
For $A$ a \DGA and $M$ a right \DG $A$-module,
the one-sided bar construction 
of \Cref{sec:bar}
is the twisted tensor product
$M \ox_{t^A} \BA$ and for a left \DG $A$-module $N$,
the two-sided bar construction $\B(M,A,N)$
is the two-sided twisted tensor product.
More generally, for $M$ an \SHM module over $A$,
the Stasheff--Halperin cochain complex defining 
their version of $\Tor_A(\kk,M)$
is computed by the twisted tensor product
$\B A \ox_t M$
for $t\: \B A \to \B (\End M) 
\to \End M$.
In Brown's setup, $C_*(B)$ is a \DGC
with comultiplication $\s \lmt \sum \s|_{[0,p]} \otimes \s|_{[p,|\s|]}$
and 
$C_*(F)$ is a $C_*(\Omega B)$-module 
via the composite 
$C_*(\Omega B) \ox C_*(F) 
\xtoo[\nabla]{} C_*(\Omega B \x F) 
\xtoo[\w_*]{ } C_*(F)
\mathrlap.
$
\eex

\section{Wolf}\label{sec:Wolf}

Joel Wolf's thesis work~\cite[Thm.~B]{wolf1977homogeneous}
proves \Cref{thm:Cartan} for $\kk$ a field.
He replaces $K$ with its maximal torus $T$ \emph{\`a la} Baum,
then deals with the noncommutativity of \eqref{eq:two-composites}
by replacing $\l_T$ with
a formality map $f\: \C(BT) \lt \H(BT)$
of Gugenheim--May type,
going in the opposite direction.

Write $\defm\rho = B(T \inc G)$.
Then under our hypotheses,
the twisted tensor product
\[{\B\C(BG)}\ox_{t^{\C(BT)} \C(\rho)}{\C(BT)}\]
computes
$\Tor_{\C(BG)}(\kk,\C BT) \iso \H(G/T)$.
Using maps of twisted tensor products induced 
by $f$ and $\l_G$ respectively,
one gets 
\[
\Tor_{\C(BG)}(\kk,\C BT) 
\lt
\Tor_{\C(BG)}(\kk,\H BT) 
\longfrom
\Tor_{\H(BG)}(\kk,\H BT), 
\]
which a generalization of \Cref{thm:EMSS-iso}
shows are linear isomorphisms.
The last Tor is computed by
the twisted tensor product
${\B\H(BG)}\ox_{t} {\H(BT)}$
with respect to $\defm t = f\C(\rho) \smash{t^{C^*(BG)}}\l_G$
and Wolf will be done if he can show that
$t$ is equal to $\H(\rho) \smash{t^{\H(BG)}}\vphantom{X^{X^x}}$
and hence gives the classical Tor.

Now $t^{\H(BT)}\B\big(f\C(\rho)\big)\l_G
=
f\C(\rho)t^{\C(BG)}\l_G$ by naturality of the tautological 
twisting cochain.
Wolf shows one can select $\l_G$ so that for $j \geq 2$,
the images of the components $\defm{\l_j} = t^{C^*(BG)}\o \l_G|_{\B_j \H(BG)}$
lie in the ideal generated by $\cupone$-products.
As $\C(\rho)$ preserves $\cupone$-products
and $f$ annihilates them,
this means $f\C(\rho)t^{\C(BG)}\l$ vanishes on $\B_{\geq 2}\H(BG)$.
Since $\H(\rho) t^{\H(BG)}$ also vanishes on $\B_{\geq 2}\H(BG)$
by definition, it remains only to check if the restrictions to
$\B_1 \H(BG) = \desusp \ol{\H(BG)}$
agree.
We can identify these with the two paths around the square
\begin{equation}\label{eq:Wolf-square}
	\begin{gathered}
	\xymatrix@R=-0.5em@C=5em{
	& \H(BT)\ar@{<-}[rd]^{f}&\\
\ol{\H(BG)} 	\ar[rd]_{\l_1 s\-}\ar[ru]^{\H(\rho)}&& \C(BT)\\
& \C(BG) \ar[ru]_{\C(\rho)}&
	}
	\end{gathered}
\end{equation}
By construction, we have identifications
$\H(\l_1 s\-) = \Hflat(\l_G) = \id_{\ol{\H(BG)}}$ 
and $\H(f) = \id_{\H(BT)}$,
so the maps in cohomology induced 
by $\H(\rho)$ and 
$f\C(\rho)\l_{\mnn 1}\susp\-$ are 
both $\H(\rho)$.
But $\H(BG)$ and $\H(BT)$ are \DGAs with zero differential,
so the maps induced in cohomology are the maps themselves, 
meaning $f\mn\C(\rho)\l_1\susp\-$ and $\H(\rho)$ are themselves equal,
completing the proof.

\bthm[Wolf]
Let $\kk$ be a field. Then \Cref{thm:Cartan}
holds additively, even without the hypothesis on $\H(BK)$.
\ethm

\section{The cobar construction}\label{sec:cobar}

The work of Husemoller--Stasheff--Moore~\cite{husemollermoorestasheff1974} 
completing the program proposed in Stasheff--Halperin 
provides a wholesale reformulation
of differential homological algebra that allows them to reprove
many of the known quasi-isomorphisms.
As this approach is homological, 
we will need to dualize several notions.
\begin{proposition}
Given a \DGC $C$, there exists a twisting cochain
$\defm{t_C}\: C \lt \defm {\W C}$
\emph{initial} in the sense that 
for any \DGA $A$,
any twisting cochain 
$t\: C \lt A$
factors uniquely through a $\DGA$ map 
$\defm{f^t}\: \W C \lt A$ 
with $t = f^t \o t_C$.
\end{proposition}
\begin{definition}\label{def:cobar}\label{def:g}
The \DGA $\W C$ is referred to as the \defd{cobar construction}
and gives the object component of a functor $\defm\W\:\Coalgs \lt \Algs$%
~{\cite[\SS1.7]{munkholm1974emss}}.
The tautological twisting cochain $t_{(-)}\: \id \lt \W$
is a natural transformation.
As an algebra, 
the cobar construction $\W C$
is the tensor algebra on the suspension
$\susp \ol C$
of the coaugmentation coideal,
and the differential is the sum of all the differentials 
$\id^{\otimes i} \ox d_{\susp \ol C} \ox \id^{\otimes j}$
and operations
$\id^{\otimes i} \otimes \,{(\susp \otimes \susp)} \D_C\desusp \otimes \id^{\otimes j}$.
There is also a dual version of the shuffle map of \Cref{def:shuffle},
a natural transformation 
\[
\defm\gamma \: \W( C_1 \ox C_2) \lt \W C_1 \ox \W C_2
\]
\nd of functors $\Coalgs \x \Coalgs \lt \Algs$.
\end{definition}
The two functors $\W \adj \B$ form an adjoint pair~{\cite[\SS1.9--10]{munkholm1974emss}}:
\[
\begin{adjunctions}
	g_t\: \W C & A \\
	f^t\:  C & \B A\mathrlap,
\end{adjunctions}
\]
linked by the twisting cochain $t\: C \lt A$.
We will have frequent recourse to the unit and counit of the adjunction $\W \adj \B$,
\[
	\defm\h\: \id \lt \B\W
\qquad\qquad\mbox{and}\qquad\qquad
	\defm\e\: \W\B \lt \id
\]
respectively. 
These are both natural quasi-isomorphisms 
and homotopy equivalences on the level of \DG modules%
~\cite[Thm.~II.4.4--5]{husemollermoorestasheff1974}%
\cite[Cor.~2.15]{munkholm1974emss}.

\bdefn
Dualizing the diagrams defining a \DG module over a \DGA $A$
gives a notion of a \defm{\DG} \defd{comodule} over a \DGA $C$.
Given a right \DG $A$-module $M$ 
and a left \DG $A$-module $N$,
the tensor product $M \ox_A N$ can 
be understood as the cokernel of the map 
$\mu_M \ox \id_N - \id_M \ox \mu_N\: M \ox A \ox N \lt M \ox N$,
and dually, given a left \DG $C$-comodule $M$
and a right \DG $C$-comodule $N$,
one can define a \defd{cotensor product} $\defm{M \square_C N}$
as the kernel of $\D_M \ox \id_N - \id_M \ox \D_N\:
M \ox N \lt
M \ox C \ox N$.
A notion of \defd{proper injective resolution} $I_\bul$
of a \DG $C$-comodule $M$ is defined dually to 
the notion of proper projective resolution
and $\defm{\mathrm{Cotor}^C(M,N)}$ is defined as 
$H_*(I_\bul \square_C N)$.
Subject to $\kk$-flatness of $C$, $H_*(C)$, $N$, and $H_*(N)$, 
there is a homological algebraic \EMSS
$\mathrm{Cotor}^{H_* (C)}(H_* M, H_* N) \implies
\mathrm{Cotor}^{ C}( M,  N)$,
which can be computed using a one-sided cobar construction
$M \ox \W C$ as a proper injective resolution of $M$.

There is also a homological Eilenberg--Moore theorem,
stating that the square of
\Cref{fig:genpull}
induces a coalgebra isomorphism $
\mathrm{Cotor}^{ C_*(B)}( C_*X,  C_*E) \iso H_*(Y)
$
(reducing to Adams's theorem 
from \Cref{rmk:Adams}
when $X \hmt {*} \hmt E$,
so that $Y \hmt \Omega B$),
and there corresponds a homological \EMSS 
$\mathrm{Cotor}^{H_* (B)}(H_* X, H_* E) \implies H_*(Y)$.
\edefn

\brmk\label{rmk:Adams}
Adams~\cite{adams1956cobar} 
introduced the cobar construction
as a model for the loop space $\Omega B$ of a simply-connected space $B$
and proved $H_*\W C_*(B) \iso H_*(\Omega B)$ 
(we dually have $\H \B\C(B) \iso \H(\Omega B)$
via Eilenberg--Moore).\footnote{\ 
	Rivera and Zeinalian recently extended 
	this to path-connected 
	spaces~\cite{riverazeinalian2018cubical,rivera2019adams}.
}
He was motivated by a desire to rephrase his 
joint result with Hilton~\cite{adamshilton1956},
which
gives a \DGA model for 		the Pontrjagin ring $H_*(\Omega B)$ in terms of the \DGC $C_*(B)$, 
in terms of a functor $\W\:\GM \lt \GM$ which would 
in principle be iterable.

Baues~\cite{baues1981double} found
a comultiplication on $\W C_*(B)$ making it a \DG Hopf algebra
and inducing the standard comultiplication on $H_*(\Omega B)$,
and allowing one iteration.\footnote{\ 
 But this goes no further: there is no suitable diagonal on $\W\W C_*(B)$.
}
Dually, he found a multiplication on $\B \C(B)$ 
making it a \DG Hopf algebra
and inducing the cup product on $\H(\Omega B)$.
This multiplication
in fact renders $\C(B)$
a \emph{homotopy Gerstenhaber algebra}
in the sense we will discuss in \Cref{sec:HGA}.
\ermk

\section{Husemoller--Moore--Stasheff}\label{sec:HSM}

Husemoller--Moore--Stasheff's proof of \Cref{thm:Cartan}
first reduces to $K = T$ a torus following Baum,
then proves the collapse of 
the homological \EMSS of \Cref{fig:Cartan}
as follows.
Taking the canonical simplicial model $K(\Z^n,1)$
for $T$, 
they note $C_*(T)$ is a commutative \DG Hopf algebra
and there are quasi-isomorphisms 
$H_*(BT) \iso H_*\B C_*(T) \lt \B C_*(T) \iso C_*(BT)$
of divided power Hopf algebras.
Moreover, they note that the cobar construction 
of both domain and codomain carry a cocommutative \DG Hopf algebra
structure $\defm{\D^*}$ preserved by this map~\cite[Prop.~IV.6.1]{husemollermoorestasheff1974}.
If $\kk$ is chosen such that $H_*(G)$ is an exterior Hopf algebra
on exterior generators $x_j$,
there are induced \DGC maps 
$C_*(BG) \simto \B C_* G \to \B\Lambda[x_j]$,
whose adjoint \DGA maps $\W C_* (BG) \lt \Lambda[x_j]$
assemble to a map 
\[
\W\B C_*(G) \isoto \W C_*(BG) \xtoo{(\D^*)^{[\rk G]}}
\big(\W C_* (BG)\big){}^{\ox \rk G} \lt 
\bigotimes \Lambda[x_j] \iso H_*(G)
\]
with adjoint $\B C_*(G) \lt \B H_*(G)$ a quasi-isomorphism.
They show the composite
\[
	H_*(BT) 
		\lt 
	C_*(BT) 
		\lt 
	C_*(BG)  
		\lt
	\B H_*(G)
\]
factors through 
\[ 
	H_*(BG) \isoto
	\bigotimes \B\Lambda[x_j] 
	\xtoo{\smash{\iter\nabla{\rk n}}}
	\B\big(\bigotimes \Lambda[x_j]\big) 
	\isoto
	\B H_*(G)
\]
because $H_*(BT) = H_* \B C_*(T)$ is commutative~%
\cite[Prop.~IV.7.2]{husemollermoorestasheff1974},
and then the commutativity of the diagram
\[
\xymatrix@R=1em{
	H_* (BT) \ar[r]\ar[dd]	&  H_*(BG)			\ar[d]\\
							& \B H_*(G)	\\
	C_*(BT) \ar[r]			& C_*(BG)\ar[u]
}
\]
gives the desired sequence of quasi-isomorphisms
\[
\xymatrix@R=1em{
\mathrm{Cotor}^{H_*(BG)}\big(\kk,H_*(BT)\big)
\ar[r]^\sim&
\mathrm{Cotor}^{\B H_*(G)}\big(\kk,H_*(BT)\big)
\\
\mathllap{.}\mathrm{Cotor}^{C_*(BG)}\big(\kk,C_*(BT)\big)
&
\mathrm{Cotor}^{C_*(BG)}\big(\kk,H_*(BT)\big)\ar[u]_\vertsim\ar[l]^(.4875)\sim
}
\]
The \DG Hopf algebra structure can be seen as a strategy 
for compilation of maps replacing the Stasheff--Halperin 
\eqref{eq:compilation}.

After all the work poured into this long paper,
there is a minor hitch in the proof of a lemma near the very end.
This lemma, however, can be substituted 
with a result of the present author,
proven for other reasons,
and an alternative definition 
due to Munkholm (see \Cref{def:psi})
of the map $\psi$ they use in defining the 
coproducts $\D^*$ on $\W C_*(BG)$ 
and $\W C_*(BT)$.\footnote{\ 
In the proof of Proposition IV.6.1,
the commutativity of the necessary diagram 
depends on 
Proposition IV.5.7,
which one can check by hand is false.
However, this diagram can be
replaced by another relevant diagram,
which actually commutes owing to the following result:
\begin{lemma}[Carlson (unpublished)]\label{thm:WDg}
	Let $A_1$ and $A_2$ be \DGAs,
	$\nabla$ as given in \Cref{def:shuffle},
	and $\psi$ as in \Cref{def:psi}.
	Then the map $\g$ of \Cref{def:g}
	is equal to the composition
	\[	
		\W(\B A_1 \ox \B A_2) \xtoo{\W\mnn\nabla}
		\W\B(A_1 \ox A_2) \xtoo{\psi}
		\W\B A_1 \ox \W\B A_2
		\mathrlap.
	\]
\end{lemma}
%
}

\bthm
Let $\kk$ be a ring such that $H_*(G)$ and $H_*(K)$
are exterior algebras on odd-degree generators.
Then there is an isomorphism of graded modules
\[
H_*(G/K) \iso \mathrm{Cotor}^{H_*(BG)}\big(\kk,H_*(BK)\big)
\mathrlap.
\]
\ethm

\section{Munkholm}\label{sec:Munkholm}%
\label{sec:compilation}


Munkholm proves \Cref{thm:Kapovitch}
by showing \eqref{eq:two-composites} is homotopy-commutative 
in a strong sense,
assuming only that $\kk$ is a principal
ideal domain and a certain extra condition in characteristic $2$.

\begin{definition}%
	[{\cite[\SS1.11]{munkholm1974emss}\cite[\SS4.1]{munkholm1978dga}}]%
	\label{rmk:bijections}\label{def:homotopy}
	Given \DGA maps $f_0,f_1\: A' \lt A$,
	a \textit{\textcolor{RoyalBlue}{\DGA}} \defd{homotopy}
	$f_0 \hmt f_1$
	is a 
	map $h\: A' \lt A$ of degree $-1$
	such that 
	\[
	\e_{\mn A} h = 0,\qquad\quad
	h\eta_{A'} = 0,
	\qquad\quad
	d(h) = f_0 - f_1,\qquad\quad
	h\mu_{\mn A'} = \mu_{\mn A}(f_0 \ox h + h \ox f_1)\
	\mathrlap.
	\]
	Given twisting cochains $t_0,t_1\: C \lt A$,
	a \defd{twisting cochain homotopy}
	$t_0 \hmt t_1$
	is a 
	map $x\: C \lt A$ of degree $0$
	such that 
	\[
	\e_{\mn A} x = \e_C,\qquad\quad
	x\h_{A} = \h_C,\qquad\quad
	\phantom{d(h) = f_0 - f_1,}\qquad
	d(x) = t_0 \cup x - x \cup t_1
	\mathrlap.
	\]
	Given \DGC maps $g_0,g_1\: C \lt C'$,
	a \textit{\textcolor{RoyalBlue}{\DGC}} \defd{homotopy}
	$g_0 \hmt g_1$
	is a 
	map $j\: C \lt C'$ of degree $-1$
	such that 
	\[
	\e_{C'} j = 0,\qquad\quad
	j\eta_{C} = 0,\qquad\quad
	d(j) = g_1 - g_0,\qquad\quad
	\D_{C'} j = (g_0 \ox j + j \ox g_1)\D_{C}
	\mathrlap.
	\]
The bijections
\[
\Algs(\W C,A)
\iso 
\Tw(C,A)
\iso 
\Coalgs(C,\B A)
\]
preserve these notions of homotopy,
as do the functors $\W$ and $\B$.
\end{definition} 



Munkholm enhances the notion of strong homotopy commutativity
by asking that the structure map $\Phi\: \B(A\ox A) \lt \B A$
satisfy up-to-homotopy unitality, commutativity, and associativity criteria.
The canonical example is that of an authentically commutative algebra.

\bex\label{thm:CGA-SHC}
If $A$ is a \CDGA, 
then $\Phi_A \ceq \B\muA\: \B(A \ox A) \lt \B A$ makes $A$ 
an \SHCA.
The cohomology ring $\H(X_\bul;\kk)$ of a simplicial set is of this type,
and will always be considered with this \SHCA structure.
If $\rho\: A \lt B$ is a map of \CDGAs,
then $\B\rho$ is an \SHCA map
per the coming \Cref{def:SHC-map}.
\eex

\bex\label{thm:H-SHC}
If $A$ is an \SHCA,
then $\H(A)$ is a \CGA,
so $\Phi_{\H(A)} \ceq \B \mu_{\H(A)}$ gives an \SHCA structure on $\H(A)$
by \Cref{thm:CGA-SHC}.
{{The cohomology ring of an}} 
{\SHC}{{-algebra will always be endowed with this}} 
{\SHC}{{-algebra structure}}.
\eex

\begin{theorem}[{\cite[Prop.~4.7]{munkholm1974emss}}]\label{thm:SHC-cochain}
	Let $X$ be a simplicial set and $\kk$ any ring.
	Then the normalized cochain algebra $\defm{\C}(X) = \C(X;\kk)$
	admits a natural \SHCA structure $\Phi_{\C(X)}$.\footnote{\ 
		This natural \SHC structure on cochains 
		is a reinterpretation of the classical 
		Eilenberg--Zilber theorem;
		only verifying the homotopy-associativity axiom requires
		much additional work. 
	}
\end{theorem}

Munkholm defines an \SHCA map 
to be a \DGC map 
making a natural square commute,
but to define the square,
we require some auxiliary concepts.
%

\begin{theorem}[{%
\cite[Prop.~IV.5.5]{husemollermoorestasheff1974}
\cite[{$k_{A_1,A_2}$, p.~21, via Prop.~2.14}]{munkholm1974emss}%
}]%
\label{def:psi}
	There exists 
	a homotopy equivalence of cochain complexes
		\[\defm \psi\: \W\B(A_1 \ox A_2) \lt \W\B A_1 \ox \W\B A_2\]
	which is a natural transformation 
	of functors $\Algs \x \Algs \lt \Algs$
	satisfying
	\[
		(\e_{\mn A_1} \ox \e_{\mn A_2}) \o \psi 
			= 
		\e_{\mn A_1 \ox A_2}\: \W\B(A_1 \ox A_2) \lt A_1 \ox A_2		
	\]
	and reducing to the identity if $A_1$ or $A_2$ is $\kk$.
\end{theorem}

\begin{definition}[{\cite[Prop.~3.3]{munkholm1974emss}}]\label{def:internal-tensor}
	Let $A_j,B_j$ be \DGAs and 
	$g_j\: \B A_j \to \B B_j$ \DGC maps for $j \in \{1,2\}$.
	The \defd{internal tensor product} 
	$\defm{g_1 \T g_2}$ is the composite
	\[
	\B(A_1 \ox A_2) 
		\xtoo{\!\!\h\!\!} 
	\B\W\B(A_1 \ox A_2)
		\xtoo{\!\B\psi\!}
	\B(\W\B A_1 \ox \W\B A_2)
		\xtoo{\!\B(\e \,\mn\W g_1 \ox \,\mn\e\mn \, \W g_2)\!}
	\B(B_1 \ox B_2)
		\mathrlap.
	\]
\end{definition}

\nd This construction 
becomes functorial on passing to the homotopy category,
or if enough of the \DGC maps $g_j$ involved are $\B f_j$ for \DGA maps $f_j$%
~\cite[Prop.~3.3(ii)]{munkholm1974emss}.
It is related to the ordinary tensor product by the shuffle map
in the sense that
$\nabla \o (g_1 \ox g_2) = (g_1 \T g_2) \o \nabla\:
\B A_1 \ox \B A_2 \lt \B(B_1 \ox \B_2)$%
~{\cite[Lem.~4.4]{franz2019homogeneous}}.

\begin{definition}\label{def:SHC-map}
	Given \SHCAs $Z$ and $A$,
a \DGC map $g\: \B Z \lt \B A$
is said to be an \defd{\SHCA map} 
if there exists a \DGC homotopy between 
	the two paths around
	the square
	\quation{\label{eq:SHC-map}
		\begin{aligned}
	\xymatrix@C=1.35em@R=2em{
		\mathllap\B(Z\ox Z)  
		\ar[d]_{g \,\,\mn\T\,\mnn g}		
		\ar[r]^(.59){\Phi_Z}
		&
		\B Z
		\ar[d]^{g}								\\
		\mathllap\B(A \ox A)
		\ar[r]_(.59){\Phi_A}			
		\ar@{}[r]_{{\vphantom{g}}}				&
		\B A\mathrlap.
	}
		\end{aligned}
}
\end{definition}

Munkholm also generalizes the construction of \eqref{eq:compilation}:
given a finite or countable list
$(g_j\: \B A_j \lt \B A)_{0 \leq j < \kappa}$ of \DGC maps,
the \defd{compilation}
\quation{\label{eq:compilation2}
	\B\big(\mn\mnn\Tensor A_j\big) 
		\xtoo[\ \ul\Tensor \, g_j\!]{} 
	\B(\mnn A^{{\otimes} \kappa}) 
		\xtoo[\iter \Phi \kappa]{}
	\B A
}
is again a \DGC map.
This operation enjoys good homotopy properties.

\bprop[{\cite[Props.~3.9(iv) \& 4.6]{munkholm1974emss}}]\label{thm:comp-pres-hmt}
The \DGC homotopy class 
of the compiled map \eqref{eq:compilation2}
depends only on the homotopy classes of the inputs $g_j$.
\eprop

\bprop[{\cite[p.~44, top]{munkholm1974emss}}]\label{thm:postcomposition-comp}
Up to \DGC homotopy, postcomposition with an \SHCA map
commutes with compilation.
\eprop

Using this enhanced notion, Munkholm is also able to 
%
verify homotopies
by checking them on generators.

\begin{theorem}[{\cite[Lem.~7.3]{munkholm1974emss}}]%
\label{thm:lambda-SHC}
	If $A$ is an \SHCA with polynomial cohomology,
	the \DGC 
	map $\l_A\: \B\H A \lt \B A$ of \eqref{eq:lambda}
	is an \SHCA map 
	when
	\bitem
	\item the characteristic of $\kk$ is not $2$, or
	\item the characteristic of $\kk$ is $2$,
	and the cup-one--squares $x_j \cupone x_j$ 
	of the chosen generators $x_j$ of $\H A$ are zero.
	\eitem
\end{theorem}
\begin{proof}
	Recall that $\l_A$ is the compilation
	of maps $\B\big(f_j\: \kk[x_j] \to A\big)$
	where $f_j(x_j) = a_j \in A$ represents $x_j \in \H(A)$.
	One checks that $\l_A \T \l_A$ is \DGC-homotopic to
	the compilation of the maps 
	$\B\big(f_j\ot \: 
	\kk[x_j] \ot \to A \ot \big)$,
	so by \Cref{thm:comp-pres-hmt},
	$\l_A$ will be an \SHCA map 
	if and only if each~$\B f_j$ is.
	To check this,
	form the associated square \eqref{eq:SHC-map},
	postcompose $t^A$
	and check if the twisting cochains
	$t' = t^A \Phi_A \B(f_j \ox f_j)$ and
	$t'' = t^A \B(f_j \mu_{\kk[x_j]})
	= f_j \mu_{\kk[x_j]}t^{\smash{\kk[x_j]{} \ot}}$
	are homotopic.
	Write $x = 1 \ox x_j$ and $y = x_j \ox 1$.
	Munkholm~\cite[Prop.~6.2]{munkholm1974emss}
	characterizes homotopy classes of twisting cochains
	$t\: \B\kk[x,y] \lt A$
	in terms of the
	classes in $\H(A)$ of three cochains:
	$t[x]$, $t[y]$, and a certain
	cocycle $c_{12}(t)$ of degree $|x|+|y|-1$.
	One easily checks $t'(x) = t''(x) = a_j = t''(y) = t'(y)$.
	If $\chara\kk \neq 2$,
	then $\H(A)$ is evenly graded, so
	$c_{12}(t') = 0 = c_{12}(t'')$ for degree reasons.
	If $\chara\kk = 2$, it turns out
	$c_{12}(t') = 2 a_j \cupone a_j = 0$
	and  $c_{12}(t'') = a_j \cupone a_j$,
	so $t'$ and $t''$ are homotopic
	if and only if the class in $H^{2|a_j|-1}(A)$ of $a_j \cupone a_j$ 
	is zero.
\end{proof}

 
\begin{theorem}[{\cite[\SS7.4]{munkholm1974emss}}]%
	\label{thm:SHC-square}
		Given \SHCAs $A$,  $X$,
		\SHCA maps $\l_{\mn X}\:$ $\B\H X \lt \B X$ and 
					$\xi\: \B A \lt \B X$,
		if $\H A$ is a polynomial ring on countably many generators,
		$\l_A$ is the \DGC map of \eqref{eq:lambda},
		and $\Hflat \l_{\mn X}$ is the identity,
		then there exists a \DGC homotopy between the two paths
		around the central square below:
\begin{equation}\label{eq:SHC-square}
	\begin{gathered}
		\xymatrix@R=-0.5em@C=2em{
		&&&\B \H(X)\ar[rrd]^{\l_X}&\\
\B\kk[x_j] \ar[r]^(.45){\B i_j}&
\B \H(A)   \ar[rrd]_{\l_A}\ar[rru]^{\B\Hflat\xi}&&&&
\B X \ar[r]^(.55){t^X} &
X\mathrlap,\\
		&&&\B A \ar[rru]_{\xi}&
		}
	\end{gathered}
\end{equation}
where $\H(A) = \kk[x_1,x_2,x_3,\ldots]$
and $i_j\: \kk[x_j] \longinc \H(A)$ are the inclusions.
\end{theorem}
\begin{proof}
First we examine the 
twisting cochain homotopy classes of the
two composites
$t\:\B\kk[x_j] \to X$ in \eqref{eq:SHC-square} for each $j$.
Each is determined~\cite[Prop.~6.2]{munkholm1974emss} 
by the class of $t[x_j]$ in $\H(X)$,
and because $\Hflat(\B\Hflat \xi) = \Hflat \xi$
and $\Hflat\l_A$ and $\Hflat\l_X$ are identity maps,
on unravelling definitions one sees they are both $(\Hflat\xi)(x_j)$.
By \Cref{rmk:bijections},
this implies
\DGC homotopies between the composites $\B\kk[x_j] \to \B X$.
By \Cref{thm:comp-pres-hmt},
it will suffice to see that $\l_X \o \B(\Hflat \xi)$
is \DGC-homotopic to the compilation of the
$\l_X \o \B(\Hflat \xi \o i_j)$
and $\xi \o \l_A$ to the compilation of the
$\xi \o \l_A \o \B i_j$.
For the former,
recall that 
$\l_X$ and $\B\Hflat\xi$ are \SHCA maps
while 
$\id_{\B \H(A)} = \B(\bigotimes i_j)$ is 
the compilation of the $\B i_j$,
and use \Cref{thm:postcomposition-comp}.
For the latter,
note $\xi$ is an \SHCA map 
and $\l_A$ was defined as 
the compilation 
$\iter{\Phi_A}\kappa \o \B(\bigotimes f_j)$
of the $\B f_j$
for $f_j\: \kk[x_j] \to A$,
so precomposing $\B i_j$ recovers $\B f_j$.
%
%
%
\end{proof}

By \Cref{thm:SHC-cochain}, 
cochain algebras are \SHCAs,
so in the situation of the Eilenberg--Moore \cref{thm:EM},
if $X \os \chi\to B \os \pi\from E$ 
have countably-generated polynomial cohomology,
one can apply \eqref{eq:lambda} three times
to obtain 
the following:
\quation{\label{eq:Munkholm-squares-complete}
	\begin{aligned}
		\xymatrix@C=2.25em@R=3.5em{
			\B\H (X)
				\ar[d]_{\l_{\mn X}}
			&	
			\B\H (B)
				\ar[d]|(.425)\hole|(.5){\l_B}|(.575)\hole		
				\ar[l]_{\B\H (\chi)}
				\ar[r]^{\B\H (\pi) }
			&
			\B\H (E) 
				\ar[d]^{\l_E}		
			\\
			\B\C(X)
			&
			\B\C(B)
				\ar[l]^{\B\C (\chi) }
				\ar[r]_{\B\C (\pi) }
			&
			\B\C(E)
			\mathrlap.
				}	
	\end{aligned}
}
If $X$ and $E$ satisfy the additional hypotheses,
one can apply \Cref{thm:lambda-SHC},
then $\l_X$ and $\l_E$ are \SHC-maps,
so applying \Cref{thm:SHC-square} twice to 
\eqref{eq:Munkholm-squares-complete},
both squares commute via \DGC homotopies.
Applying $\W$ to 
\eqref{eq:Munkholm-squares-complete}
gives a diagram $\W$\eqref{eq:Munkholm-squares-complete}
of \DGAs commutative up to \DGA homotopy.
We want to see this induces a linear isomorphism
		$\Tor_{\H (B) }(\H X,\H E) \isoto
		\Tor_{\C (B) }(\C X,\C E)$.
The \DGAs are all of the wrong form, $\W \B A$ 
but the counit $\W \B A \lt A$ of the $\W \adj \B$
adjunction is a natural \DGA quasi-isomorphism,
so by two applications of \Cref{thm:EMSS-iso},
we can instead show $\W$\eqref{eq:Munkholm-squares-complete}
induces an isomorphism
$\Tor_{\W\B\H (B) }(\W\B\H X,\W\B\H E) \isoto
		\Tor_{\W\B\C (B) }(\W\B\C X,\W\B\C E)$.

Munkholm achieves this by re-encoding homotopies.
Recall that 
a homotopy $j\: g_0 \hmt g_1\: C \lt C'$ of maps of chain complexes 
is equivalent by a single map $C \ox I \lt C'$,
where $I$ is the complex 
$\kk\{u_{[0,1]}\} \to \kk\{u_{[0]},u_{[1]}\}$
of nondegenerate chains
in the standard simplicial structure on the interval $[0,1]$.
For \DGAs,
the dual algebra $\defm{\I} = \kk\{v_0,v_1,e\}$
of normalized cochains on the simplicial interval
plays an analogous role~\cite[Thm.~5.4, pf.]{munkholm1974emss}.
For any \DGA $A$, 
the \DGA $\I \ox A$
comes with two (non-unital) \DGA maps
$\defm{\pi_j}\: \I \ox A \epi \kk{v_j} \otimes A \simto A$,
and a \DGA homotopy $h\: f_0 \hmt f_1\: A' \lt A$
induces a \DGA map
$\defm{h^P}\: A' \lt \I \ox A$
with
$f_j = \pi_j \o h^P$.%
\footnote{\ 
	To render these maps unital and augmentation-preserving,
	which is important for any further use of the adjunction,
	replace $\I \ox A$ with
	the sub-\DGA
	$PA = \kk\big\{(v_0 + v_1) \ox 1\big\} \+ \I \ox \ol A$
	and modify the $\pi_j$ in the obvious way;
	then \DGA maps $A' \to PA$ and \DGA homotopies of \DGA maps $A' \to A$
	are in bijection.
}
%
Since $\H(\I) \iso \kk$,
one sees the $\pi_j$ are
quasi-isomorphisms.

Now 
$\W$ from \eqref{eq:Munkholm-squares-complete}
induces a map of Tors in the following way~\cite[Thm.~5.4]{munkholm1974emss}.
Given \DGA maps as in \eqref{eq:Tor-DGA-functoriality-squares},
with the squares commuting up to \DGA homotopies
represented by
	$h_M^P\: A' \lt \I \ox M'$
	and $h_N^P\: A' \lt \I \ox N$.
	Then the following diagram commutes by definition:
\begin{equation}\label{eq:Tor-DGA-homotopy-squares}
\begin{aligned}
\xymatrix@C=1.90em{
  M' \ar[r]^u			&M				&\I \ox  \mnn M \ar[l]_{\pi_0}\ar[r]^(.525){\pi_1}  		&M
  \\A' \ar[u]^{\phi_{M'}} \ar[d]_{\phi_{N'}}\ar@{=}[r]		
  &A'\ar[u] \ar[d] 		\ar@{=}[r]			
  &A'\ar[u]|(.425)\hole|{h^P_M}|(.575)\hole \ar[d]|(.425)\hole|{h^P_N}|(.575)\hole \ar[r]|(.475)\hole|f|(.525)\hole
  &A\ar[u]_{\phi_M} \ar[d]^{\phi_N}
\\
  N' \ar[r]_v			&N				&\I \ox N \ar[l]^{\pi_0}\ar[r]_(.525){\pi_1}			&N\mathrlap,
}
\end{aligned}
\end{equation}
inducing a graded linear map
$
	\Tor_f(\pi_1,\pi_1) \o \Tor_{\id}(\pi_0,\pi_0)\- \o \Tor_{\id}(u,v)
	$,
	where $\Tor_{\id}(\pi_0,\pi_0)$ is an isomorphism by \Cref{thm:EMSS-iso}.
	If $u$, $f$, and $v$ are also quasi-isomorphisms,
	the composite is similarly an isomorphism.

Applying this to $\W$\eqref{eq:Munkholm-squares-complete},
we get Munkholm's main result.

\begin{theorem}
\label{thm:Munkholm}
	In the situation of \Cref{thm:EM},
	suppose that $\H(X)$, $\H(B)$, and $\H(E)$ 
	are polynomial rings on at most countably many generators,
	and if $\chara \kk = 2$, 
	assume there exist polynomial generators for $\H(X)$ and $\H(E)$
	whose $\cupone$-squares vanish.
	Then there is a graded linear isomorphism
	\[
		\Tor_{\H (B) }(\H X,\H E)
			\isoto 
		\H(Y)
		\mathrlap.
	\]
\end{theorem}	

\bcor\label{thm:Munkholm-cor}
Let $\kk$ be a principal ideal domain.
Then \Cref{thm:Cartan} holds additively.
\ecor

Munkholm observes that applying Baum's reduction \ref{thm:Baum-lemma}
and his own \cref{thm:Munkholm-cor} for $K = T$ a torus,
one recovers the Gugenheim--May \cref{thm:GuMa}.

\begin{counterexample}\label{ex:RPi}
	One might hope the isomorphism of \Cref{thm:Munkholm}
	is multiplicative, but it is not without
	some added conditions.
Let $B$ be the Eilenberg--Mac Lane space $ K(\Z/2,2)$
and $E$ the contractible space $PB$
of paths in $B$ starting at a fixed basepoint,
with $PB \to B$ evaluation at the other end.
If we take for $X$ the basepoint of $B$, 
the pullback $Y$ is the loop space 
$\Omega B = K(\Z/2,1) \hmt \Rpi$,
with polynomial cohomology ring $\F_2[\iota_1]$ over $\kk = \F_2$.
On the other hand,
$\H K(\Z/2,2)$
is a polynomial ring on generators
of degrees $2,3,5,9,17,\ldots$,
so that 
$\Tor_{\H K(\Z/2,2)}(\F_2,\F_2)$
is exterior on generators of degrees $1,2,4,8,16,\ldots$.
\Cref{thm:Munkholm} applies
and correctly reflects that the underlying 
graded $\F_2$-modules
are isomorphic,
but this isomorphism is not multiplicative.
 
\end{counterexample}

\section{Singhof}\label{sec:Singhof}

Singhof~\cite{singhof1993} directly applies Munkholm's \cref{thm:Munkholm}
to \Cref{fig:Kapovitch}
in the case
$H \x K$ acts freely on $G$, 
so that $H \lq G / K$ is a smooth manifold.

\bthm[Singhof, 1993]
Let $\kk$ be a principal ideal domain.
Then \Cref{thm:Kapovitch} holds additively.
\ethm
\bcor
If additionally $\rk G = \rk H + \rk K$, then
\Cref{thm:Kapovitch} holds multiplicatively as well.
\ecor
\bpf
In this case one has \[\Tor = \Tor^0 = \H(BH) \ox_{\H(BG)} \H(BK)\mathrlap,\]
so the ring map $\H(BH) \ox \H(BK) \to   \Tor$
is surjective.
\epf
\nd

	Singhof 
	shows the total Pontrjagin $p$ class 
	of the tangent bundle of $K \lq G /H$ 
	lies
	in the image of $\Tor^0$.
	Write
	$n = \rk G - \rk(H \x K)$
	and $n' = \dim H \lq G / K - n$.
	Through an inductive algebraic argument,
	Singhof shows 
	$\Tor^{p,q}$
	vanish outside
	the parallelogram with vertical
	edges $(p,q) \in {0} \x [0, n']$
	and ${-n} \x [-2n,n'-2n]$.
	In particular $p_\ell$ vanishes for $\ell > n'$.
	He further computes $\chi(H\lq G/K)$,
	which can nowadays be more easily done using Kapovitch's model
	in \Cref{sec:Kapovitch}.

\section{Gugenheim--Munkholm}\label{sec:GuMu}%

Gugenheim and Munkholm's joint work~\cite{gugenheimmunkholm1974}
precedes Munkholm's solo effort 
and contains arguments
that a homotopy-commutative diagram like \eqref{eq:Munkholm-squares-complete}
induces an additive isomorphism of Tors.
It also states Munkholm's hypotheses for 
the diagram to homotopy-commute, at the end,
but does not include a proof that they suffice,
and the statement that \Cref{thm:Munkholm} follows.

The main material, of which this suggested proof of \Cref{thm:Munkholm}
is an application,
is focused on extending the definition
and functoriality of Tor to a more general context than
a span $M \from A \to N$ of \DGA maps.
The approach in Munkholm to this
is that \DGC maps $\B M \from \B A \to \B N$
are enough, for applying $\W$
one gets $\W\B M \from \W\B A \to \W\B N$,
defining $\Tor_{\W\B A}(\W\B M,\W\B N)$,
and applying $\Tor_\e(\e,\e)$, 
this specializes to the $\Tor_A(M,N)$ one already had if one
was lucky enough to have \DGA maps $M \from A \to N$ to begin with.

Gugenheim--Munkholm, by contrast,
assume $M$ and $N$
are respectively right and left \DG $A$-modules
and use the two-sided bar construction $\B(M,A,N)$ as a model 
for $\Tor_A(M,N)$.
They generalize the notion of a map of \DG modules,
for $M$ a right \DG $A$-module and $M'$ a right \DG $A$-module,
and $f\: \B A \lt \B A'$ a \DGC map,
to what they call a \defm{$f$}-\defm{\textsc{sh}} \defd{linear} 
map $M \To M'$
This is a cochain map 
$g\: \B(M,A,\kk) \lt \B(M',A',\kk)$
such that $g \ox f$
makes the square with the comodule structure maps
$\B(M,A,\kk) \to \B(M,A,\kk) \ox \B A$ and
$\B(M',A',\kk) \to \B(M',A',\kk) \ox \B A'$
commute.
There is a symmetric notion 
of $f$-\textsc{sh} linear map $h\: N \to N'$
for a left \DG $A$-module $N$
and a left \DG $A'$-module $N'$,
together inducing a chain map 
$g  \ox^f h\: \B(M,A,N) \to \B(M',A',N')$
by the cotensor product and hence a map of Tors.

Gugenheim--Munkholm then introduce an appropriate notion of homotopy
of $f$-\textsc{sh} linear maps,
and show that given only homotopy-commutative squares
as in \eqref{eq:Munkholm-squares-complete} ($\l_B$ standing in for $f$),
one can actually \emph{replace} the outer maps 
by maps which are $\l_B$-\textsc{sh} linear.
Then the functorality of Tor from the previous paragraph,
will give Munkholm's \cref{thm:Munkholm},
so long as the homotopies are known to exist.

%

\section{The surjection operad}\label{sec:HGA}
The remaining work requires further cochain-level operations
generalizing
the cup-$i$ products.
The value $c'(\s|_{\D^{[0,p]}})c''(\s|_{\D^{[p,p+q]}})$
of the cup product of 
homogeneous cochains $c',c'' \in \C(X)$
and $c'' \in C^q(X;\kk)$ 
on a singular simplex $\s\: \D^{p+q} \lt X$
can be seen as the value of
$\mu_{\kk}(c' \otimes c'')$
on the 
diagonal 
$\sum_{j = 0}^{p+q} \s|_{\D^{[0,j]}} \otimes \s|_{\D^{[j,p+q]}}$;
the only term with the right dimensions
not to evaluate to zero is 
$\s|_{\D^{[0,p]}} \otimes \s|_{\D^{[p,p+q]}}$.
The higher Steenrod cup-$i$ products are defined
by apportioning the vertices differently;
for example, to define $c' \,{\cup_{\!\!\mn 1}} \,c''$
on $\s\: \D^r \to X$,
one sums, over all possible subdivisions 
$0 \leq p \leq q \leq r$,
the product of the values of $c'$ and $c''$ 
at the respective restrictions of $\s$ to 
$\D^{[0,p]\union[q,r]}$ and $\D^{[p,q]}$.\footnote{\ 
	We ignore a sign here.
}
More generally, 
given $n$ cochains $c^{(i)}$
and an $\ell$-simplex,
one can break up the vertex set $[0,\ell]$ 
into $m \geq n$ endpoint-overlapping intervals 
(in all possible ways)
and assign some of the intervals
to each $c^{(i)}$
(getting zero unless $|c^{(i)}|+1$ 
is the number of vertices assigned).
Multiplying over $i$
and summing over subdivisions then yields a cochain.
The resulting \defd{interval-cut operations}
$\C(X)^{\otimes n} \lt \C(X)$
are parameterized by the 
surjections $\varphi\:\{1,\ldots,m\} \lt \{1,\ldots,n\}$
assigning subintervals to cochains.

\begin{example}
	Identifying $\varphi$ with the sequence 
	$\big(\vp(1),\ldots,\vp(m)\big)$,
	the cup product corresponds 
	up to sign with the sequence $(1,2)$, 
	the cup-$1$ product to $(1,2,1)$,
	and the cup-$i$ product in general to an alternating sequence
	$(1,2,1,2,\ldots)$ of length $i+2$.
\end{example}

One considers only sequences with no two consecutive entries equal,
to avoid producing degenerate simplices,
and redefines the notion of cochain accordingly.
	The \defd{normalized cochain algebra}
	$\defm{\C(\simpset;\kk)}$
	on a simplicial set $\simpset$
	is the \DG subalgebra 
	containing all and only cochains 
	vanishing on each degenerate simplex.

It has been shown that the interval-cut operations
on $\C(X_\bul;\kk)$
are closed under the action of the symmetric group 
and composition,
and that the differential of such an operation
is a linear combination of other such operations%
~\cite[Prop.~1.2.7]{bergerfresse2004operad}%
\cite[Props.~2.18,\,19,\,26]{mccluresmith2003}.

In the language of operads\footnote{\ 
	Here, we use the terminology as a kind of shorthand,
	and as a familiar reference point for those familiar
	with operads,
	but for the reader who has not,
	it is enough to understand 
	``a system of composable operations (of varying arities)
	closed under composition.''
	These satisfy a long but intuitive set of axioms
	and give a formal organizing notion for certain algebraic structures.
	}	
they form a symmetric \DG-operad~$\defm{\ms X\!\mnn}$, 
called the 
\defd{surjection operad}~\cite{bergerfresse2004operad}.
	The normalized cochain algebra
	of a pointed simplicial set is then functorially
	a $\ms X\!\mnn$-algebra~{\cite[Thm.~2.15]{mccluresmith2003}}.

In 1963, Gerstenhaber~\cite{gerstenhaber1963cohomology}
observed that  
the desuspension of the Hochschild cohomology ring of an associative
algebra carried
an additional graded Lie algebra structure 
rendering the operations $[a,-]$ (for homogeneous $a$)
graded derivations;
such a structure became known as a \emph{Gerstenhaber algebra}.
In the early 1990s, 
it was observed by various authors%
~\cite{getzlerjones,gerstenhabervoronov1995},
including Gerstenhaber himself,
that the 
operations on the Hochschild complex
used to define the Gerstenhaber bracket
extended to a more complicated
algebraic structure which 
would be
connected to the surjection operad.

\begin{definition}[{\cite[\SS3.2]{franz2019shc}}; %
	{\cite[\SS2.1]{kadeishvili2003cochain}}	over $\F_2$;
	{\cite[(3.13)]{franz2020szczarba}} for signs]
	\label{def:xHGA}\label{def:HGA-map}
\ 	\\ 
A \defd{homotopy Gerstenhaber algebra} 
	(\textcolor{RoyalBlue}{\HGA})
	is an algebra over the symmetric \DG operad
	$F_2 \ms X\!\mnn$ of~$\ms X\!$
	generated by the \DGA product $(1,2)$ 
	and the operations
	$\defm{E_\ell}$ 
	corresponding to $(1,2,1,3,1,\ldots,1,\ell+1,1)$.%
	\footnote{\
	\HGAs were first defined~\cite{gerstenhabervoronov1995}
		as \DGAs equipped with operations $E_\ell$
		satisfying certain axioms.
		It is a theorem%
			~\cite[Thm.~4.1]{mccluresmith2003}%
			\cite[\SS1.6.6]{bergerfresse2004operad}
		that this yields precisely $F_2 \ms X\!$-algebras.
	}
	An \defd{extended} \defm{\HGA}  
	is an algebra over the symmetric \DG operad
	 $\defm{F'_3 \ms X}\!$ of $\ms X\!$ 
 	generated by $F_2 \ms X\!$ and the
	operations $\defm{F_{p,q}}$
	corresponding to 
	\[
		(
		1,\ p + 1,\,
		1,\ p + 2,\,
		1,\ p + 3,\,
		\ldots,\,
		1,\ p + q,\,
		\, 
		1,\ p + q,\,
		2,\ p + q,\,
		3,\ p + q,\,
		\ldots,\,
		p,\ p + q
		)\mathrlap.
	\]
	An (\defd{extended})
		\textcolor{RoyalBlue}{\HGA} 
		\defd{homomorphism} 
	$f\: A \lt B$ 
	is a \DGA map
	distributing over the operations ($F_{p,q}$ and) $E_\ell$.
\end{definition}

To explain the notation, the \DG suboperad $F_2 \ms X$ 
is a term of a certain increasing filtration $F_n \ms X$
on $\ms X$ with $F'_3 \ms X \sub F_3\ms X$.%
\footnote{\ 
	$F_0 \ms X\!$ contains functions with $0$ or $1$ value.
	$F_\ell \ms X \less F_{\ell-1} \ms X\!\mnn$
	contains sequences such that $\ell$ is 
	the maximum number of alternations 
	in a two-value subsequence.
	For example, $(1,2)$ lies in $F_1\ms X\!$
	and $(3,1,4,3,4,2,1,2)$ in $F_3\ms X\!$,
	with maximum alternation 
	attained by the subsequence $(3,1,3,1)$.
	}$^,$\footnote{\ 
	For those who know what this is, $F_n \ms X$ is equivalent to the cochains 
	on the operad of little $n$-cubes.
	In fact~\cite[Thm.~1.3.2, Lem.~1.6.1]{bergerfresse2004operad},
	$\ms X\!$ is a quotient of the 
	\DG-operad $\defm{\ms E}$ associated
	to the classical Barratt--Eccles simplicial operad,
	which is filtered by a sequence $F_n \ms E$ 
	of $E_n$-operads
	(of which the $F_n \ms X\!$ are the images),
	so the normalized cochain algebra is an $E_\infty$-algebra.

	The filtrands $F_n \ms X\!$ had been
	already identified and shown to be equivalent
	to the little $n$-cubes operads by other methods~\cite{mccluresmith2003}
	before this surjection
	was found. 
	Earlier still,
	McClure--Smith 
	had shown
	$F_2 \ms X\!$ is equivalent to the little squares 
	operad~\cite{mccluresmith2002deligne}
	in order to prove Deligne's conjecture
	that the Hochschild cohomology of a ring 
	is naturally an $E_2$-algebra.
}

\begin{example}\label{ex:CDGA-surjection}
	A \CDGA $A$ is canonically an $\ms X\!\mnn$-algebra
	with $F_\ell \ms X\less F_{\ell-1} \ms X\!$ 
	acting identically
	as $0$ for $\ell \geq 2$.
	In particular,
	a \CDGA $A$ is naturally an extended \HGA.
\end{example}

\begin{corollary}\label{thm:cochain-xHGA}
		For any pointed simplicial set $X_\bul$, 
		its algebra $\C(X_\bul)$
		of normalized cochains is also naturally an extended \HGA
		with $-E_1 = \cupone$ and $-F_{1,1} = \cuptwo$.
\end{corollary}

	An \HGA structure on a \DGA $A$ is known to induce 
	a multiplication $\mu_{\B A}$ on the \DGC $\B A$
	rendering it a \DG Hopf algebra.\footnote{\
		Extended \HGAs over $\F_2$
		were first studied by 
		Kadeishvili~\cite{kadeishvili2003cochain},
		seeking conditions on an \HGA $A$ 
		under which $\B A$ admitted cup-$i$ products.
		He gave a characterization in terms of 
		operations of $\ms X\!$ acting on $A$;
		for $i = 1$, these are the $F_{p,q}$.
	}

\begin{theorem}[Franz~{\cite{franz2019shc}}]\label{thm:Franz-SHC}
A homotopy Gerstenhaber algebra $A$ admits a 
natural \DGC map \[\Phi_A\: \B(A \ox A) \lt \B A\]
satisfying the unitality and associativity axioms for \SHCAs
and such that $\mu_{\B A} = \Phi_A \o \nabla$.
If $A$ is an extended homotopy Gerstenhaber algebra,
then $\Phi_A$ is an \SHCA structure.
\end{theorem}
 

This \SHCA structure is not complicated to define,
but the fact it actually satisfies the axioms
is an extraordinarily complex computational result 
(found by hand but verified by computer)
for which we still have no conceptual explanation.

\section{Franz}\label{sec:Franz}

In 2019,
Franz~\cite{franz2019homogeneous} 
proved the following.
\bthm[Franz]
Let $\kk$ be principal ideal domain in which $2$ is a unit.
Then \Cref{thm:Cartan} holds.
\ethm

Recall that the characteristic-$0$ multiplicative
results work because there exist \DGA structures on the complexes
computing $\Tor$.
Kadeishvili--Saneblidze%
~\cite[Thm.~7.1]{kadeishvilisaneblidze2005cubical}
found a \DGA structure on
the twisted tensor product $\B A \ox_t A'$ 
when $A$ and $A'$ are \HGAs
and $t = t^{A'} \o \B\phi$ for an \HGA homomorphism
$\phi\: A \to A'$.
This \DGA structure is functorial on the category 
of \HGA diagrams of  shape $\bul \to \bul$.
By 
\Cref{thm:EM},
\Cref{def:one-sided}, and
\Cref{thm:cochain-xHGA},
$\B \C(BG) \ox_{t^{\C(BK)}\B\rho} \C(BK)$
is a \DGA with cohomology ring $\H(G/K)$,
and
$\B \H(BG) \ox_{t^{\H(BK)}\B\H\rho} \H(BK)$
is a \DGA with cohomology ring
$\Tor_{\H BG}(\kk,\H BK)$.
Using the Halperin--Stasheff maps \eqref{eq:lambda},
Franz considers the composite cochain map
\[
\xymatrix{
	\defm{\Theta}\:
\dsp\B \H(BG)\quad \ox_{{\mathclap{t^{\H(BK)}\,\B\H\rho}}}	\quad\H(BK)
	\ar[r]^(.539){\id \ox \lho}&
\dsp\B \H(BG)\quad \ox_{\mathclap{t^{\C(BK)}\,\l_H\,\B\H\rho}} 	\quad\C(BK)
	\ar[d]_(.55){\vertsim}^(.57){\d_h}\\
\dsp\B \C(BG)\quad \ox_{\mathclap{t^{\C(BK)}\,\B\rho}} 			\quad\C(BK)	
	&
\dsp\B \H(BG)\quad \ox_{\mathclap{t^{\C(BK)}\,\B\rho\,\l_G}} 	\quad\C(BK)
	\ar[l]^(.495){\l_G \ox \id}
}
\]
where $\defm\lho = t^{\H(BK)} \o \l_H \o \desusp_{\C(BK)}$
and $\defm{\d_h}$ is the cap product with 
the homotopy~$h$ of twisting cochains
$\B\C(BG) \lt \H(BK)$ implied by the commutativity of 
\eqref{eq:two-composites},
as proved by Munkholm.
By a variant of \Cref{thm:EMSS-iso},
$\H(\Theta)$ 
is a linear isomorphism.

Because $\l_G$ and $\l_H$ are not multiplicative,
$\Theta$ may not preserve the \DGA structures,
but we follow Wolf
in postcomposing a formality map
to kill $\cupone$-products
preventing $\l_G$ from being multiplicative.
Let $T$ be a maximal torus of $H$ 
and set $\phi = B(T \inc H)$.
Franz follows $\Theta$ with 
\[
\defm{\Psi}\:
\B \C(BG)  \ox_{{t^{\C(BK)}\,\B\rho^*}} 	\C(BK)	
	\xtoo{\id \ox f\mn \phi^*}
\B \C(BG) \ox_{{t^{\C(BT)}\,\B(f\phi^*\rho^*)}} 	
\H(BT) \mathrlap,
\]
where
$\id \ox \phi^*$ 
induces an injection in cohomology
and $\id \ox f$
is a quasi-isomorphism by \Cref{thm:EMSS-iso}.
Since $\phi$ is an \HGA map,
if $f$ can also be chosen to be an \HGA map,
$\Psi$ will be multiplicative
by functoriality of the \DGA structure.

The existing Gugenheim--May map 
does not 
meet these desiderata, so Franz 
defines a certain ideal \mathversion{normal2}$\defm{\f k_X}$ 
of undesirable cochains,
functorial in spaces~$X$,
and then constructs a new formality map $\defm f\: \C(BT) \lt \H(BT)$ 
which is an \HGA map annhilating $\f k_{BT}$
\emph{so long as $2$ is a unit of $\kk$}.
\mathversion{normal}
At the same time, using the fact $\C(BT)$ is an extended \HGA,
he applies \Cref{thm:Franz-SHC}
to obtain an \SHCA structure $\Phi$
such that
the \DGC quasi-isomorphism $\l_T \: \B\H(BT) \lt \B\C(BT)$ 
compiled with respect to $\Phi$ has image
in the error ideal\mathversion{normal2} $\f k_{BT}$.
Moreover, the compiled maps make \eqref{eq:two-composites}
commute up to a homotopy $h$ 
such that $t^{\C(BK)}h$ 
is congruent to the identity 
modulo $\f k_{BK}$.
Further, $\l_G$ is an \SHCA map
in the sense of making \eqref{eq:SHC-map} 
homotopy-commute
via a \DGC homotopy
whose associated twisting cochain homotopy 
takes the coaugmentation coideal
into $\f k_{BG}$.
This twisting cochain homotopy
is hence annihilated by postcomposing with 
$f\mn\phi^*\rho^*$.\mathversion{normal}
With these choices of \SHCA structure and formality map, 
one checks postcomposing with $\Psi$
simplifies $\Theta$
so that the composite $\Psi\Theta$ is
\[
\l_G \ox \phi^*\: 
\B\H(BG) \ox_{\H(\rho) t^{\H(BG)}} \H(BK) 
	\lt 
\B\C(BG) \ox_{f\mn\phi^*\rho^*t^{\C(BG)}} \H(BT)
\mathrlap.
\]

Since $\H(\Psi)$ is a multiplicative injection and $\H(\Theta)$
is a bijection, to show  $\H(\Theta)$ is multiplicative,
it suffices to show $\H(\Psi\Theta)$ is.
When $A'$ is a \CDGA, 
the multiplication on the 
twisted tensor product $\B A \ox_t A'$
making it a \DGA 
is just the naive multiplication 
permuting the tensor factors and then multiplying components
using $\mu_{\B A} \ox \mu_{A'}$,
because the commutativity of $A'$
means the higher-order \HGA operations on $A'$ 
figuring in the multiplication formula vanish.
To show $\H(\Psi\Theta)$ is multiplicative,
it is thus enough to construct a \DGC homotopy between
the two paths around the large diagram
\[
\resizebox{125mm}{!}{%
\xymatrix{
	(\B\H BG)\ot \otimes (\H BK) \ot \ar[r]^(.5){\nabla \ox \id}\ar[d]_{\l_G\ot \,\ox\  (\phi^*)\ot}&
	\B(\H BG)\ot \otimes (\H BK)\ot \ar[r]^(.555){\Phi \ox \mu}\ar[d]_{\l_G\T\l_G \,\ox\  (\phi^*)\ot}
	&\B\H BG  \otimes \H BK  \ar[d]^{\l_G \ox \phi^*}\\
	(\B\C BG)\ot \otimes (\H BT)\ot \ar[r]_(.5){\nabla \ox \id}
	&
	\B(\C BG)\ot \otimes (\H BT)\ot \ar[r]_(.555){\Phi \ox \mu}
	&\B\C BG  \otimes \H BT \mathrlap.
	}
}
\]
By \Cref{thm:Franz-SHC}, the product on $\B\C(BG)$
factors as $\Phi \o \nabla$,
so it is enough to find a \DGC homotopy 
for the right square.
The \DGA coordinates commute on the nose
since $\phi^*$ is a ring map.
For the bar coordinates,
$\l_G$ is an \SHCA map
via a homotopy $H$ such that 
the coaugmentation coideal of $\B (\H BG)\ot$
is annihilated by
the twisting cochain $t^{\H(BT)}f\mn \phi^*\rho^*H$ 
defining the twisted differential on the codomain.
An easy lemma then shows 
$H \ox \phi^*$ is the desired \DGC homotopy.

\begin{example}\label{eg:recurrent}
	Let $H \iso \U(1)$ be the subgroup of $\SU(4)$ with diagonal entries ${\diag}(z^{-3},z,z,z)$.
	Then,
	indexing generators by degree, 
	we have a ring isomorphism
	\[
	\H\big(\SU(4)/H;\Z[{\textstyle\frac 1 2}]\big) \,\iso\, \frac{\Z[\frac 1 2][s_2] \otimes \Lambda[a_5,b_7]}
	{(3s^2,s^3,s^2 a_5)}\mathrlap.
	\]
\end{example}

\section{Carlson--Franz}\label{sec:CF}\label{sec:two-sides}

On seeing Franz's paper,
the present author at once had the insight that Franz's proof 
of \Cref{thm:Cartan}
would effortlessly generalize
to a proof of \Cref{thm:Kapovitch}.
That insight was wrong:
Franz pointed out no one had defined a multiplication $\defm\mu$ on 
a two-sided twisted tensor product $A'' \ox_{t''} \B A \ox_{t'} A'$
generalizing the product of Kadeishvili--Saneblidze.
One also 
needs a version of the Eilenberg--Moore theorem preserving $\mu$;
more precisely,
one needs to show that in the situation of \Cref{thm:EM},
the relevant instance
$
	\C(X) \ox_{t''} \B\C(B) \ox_{t'} \C(E) 
		\to 
	\C(X) \ox \C(E) 
		\to
	\C(Y)
$
of the map of \eqref{eq:EM-res}
is multiplicative up to a homotopy $\defm h$.

Within a month, however, 
Franz shared formulas for $\mu$ 
and $h$
he believed should work
once a correct choice of signs was discovered.
The present author guessed the signs
and proved that these formulas 
work as expected.\footnote{\ 
	After producing the adapted 
	variant of the Eilenberg--Moore theorem,
	one also needs to extend some lemmas on twisted 
	tensor products to
	the two-sided case, 
	and to check that a map defined by Wolf
	does what is required 
	(this is harder,
	and proof formerly appeared only in Wolf's unpublished thesis)%
	---but 
	the product on the two-sided bar construction
	is the important idea.
	}
After that, 
the proof summarized in the preceding section 
works \emph{mutatis mutandis}.
This was the work presented at the conference session.

\bthm[{\cite{carlsonfranzlong}}]\label{thm:CF}
\Cref{thm:Kapovitch} holds whenever $2$ is a unit of $\kk$.
\ethm

\nd As one needs to invert $2$ for this approach to work,
this result is not strictly an improvement on Munkholm's
additive \cref{thm:Munkholm} in all cases.

\bs 

\nd\emph{Acknowledgments.} 
This work was conducted during a visiting scholarship at Tufts,
which the author would like to gratefully acknowledge.
The author would further like to thank
Markus Szymik for a copy of 
Stasheff--Halperin~\cite{halperinstasheff1970},
Paul Baum, Peter May, Larry Smith, Jim Stasheff, and Joel Wolf
for first-hand historical insight,
Jim Stasheff and the referee for writing advice,
and Joanne Quigley for proofreading.

 	{\footnotesize\bibliography{bibshort}}

	\bigskip

\nd\footnotesize{\url{jeffrey.carlson@tufts.edu}}
\end{document}